\documentclass[sn-mathphys-num,iicol]{sn-jnl}
\usepackage{preprint-notation}

\usepackage{graphicx}%
\usepackage{multirow}%
\usepackage{amsmath,amssymb,amsfonts}%
\usepackage{stmaryrd}%
\usepackage{amsthm}%
\usepackage{mathrsfs}%
\usepackage[title]{appendix}%
\usepackage{xcolor}%
\usepackage{textcomp}%
\usepackage{manyfoot}%
\usepackage{booktabs}%
\usepackage{algorithm}%
\usepackage{algorithmicx}%
\usepackage{algpseudocode}%
\usepackage{listings}%
\usepackage{tikz}%
\usetikzlibrary{arrows.meta,decorations.pathreplacing,calc}

\theoremstyle{thmstyleone}
\newtheorem{lemma}{Lemma}
\newtheorem{proposition}[lemma]{Proposition}
\newtheorem{corollary}[lemma]{Corollary}
\newtheorem{assumption}[lemma]{Assumption}
\theoremstyle{remark}
\newtheorem{rem}{Remark}

\begin{document}
\title{Skeleton-Stabilized CutIGA for the Stokes Problem} 

\author*[1,2]{\fnm{E. Harald} \sur{van Brummelen}}\email{e.h.v.brummelen@tue.nl}
\author[3]{\fnm{Mats G.} \sur{Larson}}\email{mats.larson@umu.se}
\author[1]{\fnm{Clemens V.} \sur{Verhoosel}}\email{c.v.verhoosel@tue.nl}
\author[1,4]{\fnm{Hugo M.} \sur{Verhelst}}\email{h.m.verhelst@tue.nl}

\affil*[1]{\orgdiv{Department of Mechanical Engineering},
  \orgname{Eindhoven University of Technology},
  \orgaddress{\street{P.O. Box 513}, \postcode{5600 MB}, \city{Eindhoven},
  \country{The Netherlands}}}

\affil[2]{\orgdiv{Department of Mathematics and Computer Science},
  \orgname{Eindhoven University of Technology},
  \orgaddress{\street{P.O. Box 513}, \postcode{5600 MB}, \city{Eindhoven},
  \country{The Netherlands}}}

\affil[3]{\orgdiv{Department of Mathematics and Mathematical Statistics},
  \orgname{Ume\aa{} University},
  \orgaddress{\postcode{SE-901 87}, \city{Ume\aa}, \country{Sweden}}}

\affil[4]{\orgdiv{Faculty of Military Sciences},
  \orgname{Netherlands Defence Academy},
  \orgaddress{\street{Het Nieuwe Diep 8}, \postcode{1781 AC}, \city{Den Helder},
  \country{The Netherlands}}}

\date{\today}

\abstract{%
\unboldmath%
We present a stability, error and conditioning analysis of skeleton-stabilized
cut isogeometric analysis (CutIGA) for the Stokes problem. The physical domain
is embedded in a structured background mesh, and the velocity and pressure are
approximated by equal-order spline spaces, with maximum regularity as the
primary case and reduced regularity down to standard $C^0$ elements included.
The inf-sup deficiency of the equal-order pairing is compensated by a skeleton
stabilization that penalizes jumps in high-order normal derivatives of the
pressure across all interior faces of the active mesh, and robustness with respect to arbitrarily
small cut elements is recovered by a ghost penalty for the velocity on the
faces of cut elements. Dirichlet conditions are imposed weakly by Nitsche's
method, and Neumann boundaries may intersect the mesh arbitrarily. For
sufficiently fine meshes, we prove inf-sup stability, optimal-order a priori
error estimates, and an $O(h^{-2})$ bound on the condition number of
the system matrix, with constants independent of the cut configuration.
Numerical experiments for a channel with a cylinder, rotated with respect to
the background mesh, confirm the predicted convergence rates and the
cut-independent conditioning.}


\maketitle

\section{Introduction}
\label{sec:intro}

Isogeometric analysis (IGA) was introduced by Hughes \emph{et al.}~\cite{Hughes2005} with the objective of integrating Computer-Aided Design (CAD) and finite-element analysis. In IGA, the spline functions used to describe the computational geometry are also employed to construct the approximation spaces for the governing partial differential equations. Besides enabling design-through-analysis workflows, spline-based discretizations possess several attractive properties, including higher inter-element continuity, superior approximation properties, and improved approximation of higher-order derivatives~\cite{Cottrell:2009fk}. Among the many combinations of polynomial degree and regularity offered by spline spaces, maximum-regularity spline spaces are particularly attractive owing to their superior approximation efficiency on a per-degree-of-freedom basis~\cite{Evans:2009kq}. 

Despite considerable progress in geometric modeling and computational analysis, fully automated procedures for generating boundary-conforming volumetric analysis models directly from CAD data remain a topic of ongoing research for general geometries. Unfitted methods, including CutFEM~\cite{Burman:2025sf,Burman:2015fp}, CutIGA~\cite{Elfverson2018}, the finite-cell method~\cite{Rank:2012ye}, and immersogeometric analysis~\cite{Kamensky:2015zh,Hoang:2019ty}, circumvent this bottleneck by embedding the physical domain in a geometrically simple ambient domain discretized with a structured mesh. The approximation space is constructed on the background mesh and subsequently restricted to the physical domain, thereby retaining the approximation properties of the underlying spline discretization while eliminating the need for boundary-conforming meshes. 
This approach, however, introduces several challenges, including accurate numerical integration over cut elements~\cite{Loibl:2026ys,Burman:2025sf}, the weak imposition of boundary conditions~\cite{Nitsche:1971fk,Ruess:2013fk,Burman:2012nn}, and the loss of stability associated with arbitrarily small cut elements~\cite{Prenter:2023rp,Prenter:2017rw}.

Besides the aforementioned geometric considerations, 
approximation methods for mixed problems such as the Stokes equations must account for stability requirements associated with the pairing of approximation spaces. Stable discretizations of the Stokes problem are traditionally obtained by pairing velocity and pressure spaces that satisfy the discrete inf-sup (LBB) condition~\cite{BoffiBrezziFortin2013}. Within isogeometric analysis, several inf-sup stable spline discretizations have been developed, including both Taylor--Hood-type pairings and compatible spline complexes~\cite{Buffa:2011yc,Evans:2012kq}. However, inf-sup stability is generally linked to the correspondence between basis functions and the topology of the underlying mesh and, consequently, velocity-pressure pairs that are inf-sup stable on fitted meshes need not remain stable in an unfitted setting~\cite{Hoang:2019ty,GuzmanOlshanskii2018,Kirchhart2016}. Furthermore, inf-sup stable spline pairings generally do not permit maximum-regularity spline spaces for all components~\cite{Bressan:2013fy,Evans:2012kq}.

{\em Stabilization\/} provides a general strategy for resolving the diverse instabilities that may arise in finite-element formulations, including those induced by pathologically cut elements in unfitted discretizations and those arising from inf-sup unstable  approximation-space pairings in mixed formulations. Various stabilization techniques have been developed to address both classes of instabilities.

Stabilization techniques for inf-sup unstable mixed formulations have been studied extensively since the pioneering work of Brezzi and Pitk\"aranta~\cite{Brezzi:1984zr}. Their common objective is to compensate for the {\em LBB deficiency\/}~\cite{Bochev:2006gf} of the chosen approximation-space pairing by augmenting the discrete formulation with suitably constructed stabilization terms. According to the support of these terms, stabilization methods are generally divided into element-interior and edge-based~\cite{Douglas:1989hi,Burman:2006wf,Bochev:2006gf} methods. Element-interior stabilization methods can be further divided into residual-based~\cite{Hughes:1986uz,Hughes:1987uo,Douglas:1989hi,Brezzi:1988il,Tezduyar:1991by,Barth:2004gf} and projection-based~\cite{Becker:2001pb,Braack:2006wl,Bochev:2006gf} formulations. 

Stabilization of unfitted discretizations aims to eliminate the loss of robustness of the discrete formulation associated with arbitrarily small cut elements~\cite{Prenter:2017rw,Prenter:2023rp}. Such stabilization methods can generally be divided into extension-based methods and augmentation-based methods. In extension-based methods, such as the finite-cell method~\cite{Parvisian:2007eu,Rank:2012ye}, the functionals of the original weak formulation are extended onto the ambient domain. Typically, the extension employs modified parameters, such as a fictitious-domain stiffness, which introduces a consistency error that is controlled by these parameters. Augmentation-based methods such as the ghost-penalty method~\cite{Burman:2010km} instead augment the original weak formulation with weakly consistent stabilization terms, which vanish at an optimal rate for sufficiently smooth solutions, to recover robustness. Depending on the support of these stabilization terms, augmentation-based methods can be further divided into edge-based~\cite{Burman:2012nn,Burman:2015fp} and patch-based~\cite{Burman:2026ls,Burman2023ext} methods. An alternative to stabilization is provided by cell-agglomeration techniques, which recover robustness by modifying the approximation space rather than the weak formulation~\cite{Badia:2018oe}.

The stabilization techniques developed for mixed and unfitted finite-element discretizations extend naturally to isogeometric analysis. However, by virtue of the increased inter-element continuity of spline spaces, the stabilization terms in IGA can be restricted to 
higher-order derivatives. This principle underlies the skeleton-stabilization framework for mixed spline discretizations~\cite{Hoang:2018lp} and its extension to immersogeometric analysis for the Stokes problem~\cite{Hoang:2019ty}, in which skeleton stabilization is employed to compensate for the inf-sup instability of equal-order, maximum-regularity spline spaces and ghost-penalty stabilization is used to recover robustness with respect to arbitrarily small cut elements. Although numerical experiments have demonstrated the effectiveness of this formulation, a rigorous mathematical analysis establishing its stability, convergence and conditioning is presently unavailable. 

The objective of the present work is to provide the first rigorous analysis
of skeleton-stabilized CutIGA for the Stokes problem. For sufficiently fine
meshes, we establish inf-sup stability, optimal-order a priori error
estimates, and a cut-independent bound on the condition number. We focus on
maximum-regularity equal-order approximations of the velocity and pressure,
but also consider reduced-regularity approximations, down to standard $C^0$
splines (FEM). The structure of the analysis is in line with that of stabilized
Nitsche fictitious-domain methods for equal-order $P_1$--$P_1$
elements~\cite{BurmanHansbo2014,Massing:2014}; see also~\cite{Hansbo2014}.
The extension to splines of higher degree and regularity is, however, not
immediate: the skeleton stabilization only controls jumps of higher-order
normal derivatives, which we relate to a gradient-jump seminorm
in Lemma~\ref{lem:seminormequiv}, and the extended supports of the spline basis
functions require a quasi-interpolant that involves only the active domain
(Appendix~\ref{app:cutiga-tools}). In addition, Neumann boundaries that
intersect the mesh arbitrarily require a dedicated treatment
(Appendix~\ref{app:neumann-boundary-neighborhood}). For
boundary-conforming meshes, in which case the ghost penalty is not required,
the analysis also covers the skeleton-stabilized formulation
of~\cite{Hoang:2018lp}, in the variant with weakly imposed Dirichlet
conditions.

The remainder of this paper is organized as follows.
Section~\ref{sec:cutiga} introduces the Stokes problem and its skeleton-stabilized
CutIGA approximation. Sections~\ref{sec:stability}--\ref{sec:condition-number}
establish stability, a priori error estimates, and a cut-independent bound on
the condition number, respectively. Section~\ref{sec:numerical-experiments}
presents numerical experiments that illustrate these results, and
Section~\ref{sec:conclusions} concludes. The appendices
collect the auxiliary estimates.

\section{CutIGA approximation of the Stokes problem}
\label{sec:cutiga}
In this section we present the Stokes problem and its approximation by means of the skeleton-stabilized immersed IGA formulation in~\cite{Hoang:2019ty}.
\subsection{Stokes problem}
\label{sec:StokesProblem}
We consider a bounded Lipschitz domain $\Omega\subset\IR^d$
($d=2,3$) with boundary $\Gamma:=\partial\Omega$. The boundary admits a decomposition into complementary subsets $\Gamma_{\textsc{D}}$ and $\Gamma_{\textsc{N}}$ such that 
$\operatorname{meas}(\Gamma_{\textsc{D}})>0$ and
$\operatorname{meas}(\Gamma_{\textsc{N}})>0$; the pure Dirichlet case
$\Gamma_{\textsc{N}}=\emptyset$ is addressed in
Corollary~\ref{cor:pure-dirichlet} and Remark~\ref{rem:pure-dirichlet}. Denoting the exterior unit normal vector on~$\Gamma$ by~$\bfn$, the Stokes equations for the velocity field $\bfu:\Omega\to\IR^d$ and pressure field $p:\Omega\to\IR$ read
\begin{subequations}\label{eq:stokes}
\begin{alignat}{2}
-\nabla\cdot(2\mu\nabla_{\mathrm{s}}\bfu)+\nabla{}p&=\bff
&\quad&\text{in }\Omega
\label{eq:stokes-a}
\\
\nabla \cdot \bfu &= 0  &\quad&\text{in }\Omega
\label{eq:stokes-b}
\\
\bfu&=\bfg  &\quad&\text{on }\Gamma_{\textsc{D}}
\label{eq:stokes-c}
\\
2\mu(\nabla_{\mathrm{s}}\bfu)\bfn-p\bfn&=\bfh&\quad&\text{on }\Gamma_{\textsc{N}},
\label{eq:stokes-d}
\end{alignat}
\end{subequations}
with $\nabla_{\mathrm{s}}\bfu:=\frac{1}{2}(\nabla\bfu+(\nabla\bfu)^T)$. Here,
$\mu>0$ denotes the dynamic viscosity. To condense the presentation, we
set $2\mu=1$ throughout. The dependence on a general constant coefficient
can be recovered by scaling the viscous terms accordingly.

To provide a basis for the immersed formulation, we consider the weak formulation of~\eqref{eq:stokes} according to:
\begin{multline}\label{eq:stokes-weak-formulation}
\text{Find }(\bfu,p)\in\bfH{}^1_{\Gamma_{\textsc{D}},\bfg}(\Omega)\times
L^2(\Omega):\\
a(\bfu,\bfv)+b(\bfv,p)-b(\bfu,q)=\ell(\bfv)\\
\forall(\bfv,q)\in\bfH{}^1_{\Gamma_{\textsc{D}},0}(\Omega)\times
L^2(\Omega),
\end{multline}
with $\bfH{}^1_{\Gamma_{\textsc{D}},\bfg}(\Omega)$ the space of $\IR^d$-valued functions with components in~$H^1(\Omega)$ and trace $\bfg$ on~$\Gamma_{\textsc{D}}$ and
\begin{subequations}
\begin{align}
a(\bfv,\bfw)&:=(\nabla_{\mathrm{s}}\bfv,\nabla_{\mathrm{s}}\bfw)_{\Omega}\\
b(\bfv,q)&:=-(\nabla\cdot\bfv,q)_{\Omega}
\label{eq:stokes-bilinear-b}\\
\ell(\bfv)&:=(\bfv,\bff)_{\Omega}+(\bfv,\bfh)_{\Gamma_{\textsc{N}}},
\end{align}    
\end{subequations}
where $(\cdot,\cdot)_{\omega}$ denotes the $L^2(\omega)$ inner product for scalar-, vector- or tensor-valued functions.

\subsection{CutIGA approximation}
\label{sec:stabilized-cutiga}

To provide a setting for the CutIGA approximation of the Stokes problem, we embed $\Omega$ into an ambient box
$\widetilde{\Omega} = \prod_{i=1}^d (a_i, b_i) \supset \overline{\Omega}$, where $\overline{\,\cdot\,}$ denotes the closure. We denote by
$\mathcal{H}\subset(0,h_0]$ a countable set of mesh parameters with $0$ as its only accumulation point, and by
$\{\widetilde{\mcT}_h\}_{h\in\mathcal{H}}$ a corresponding family of
quasi-uniform, shape-regular tensor-product meshes on
$\widetilde{\Omega}$.

\begin{figure}
\centering
\includegraphics[width=0.95\columnwidth]{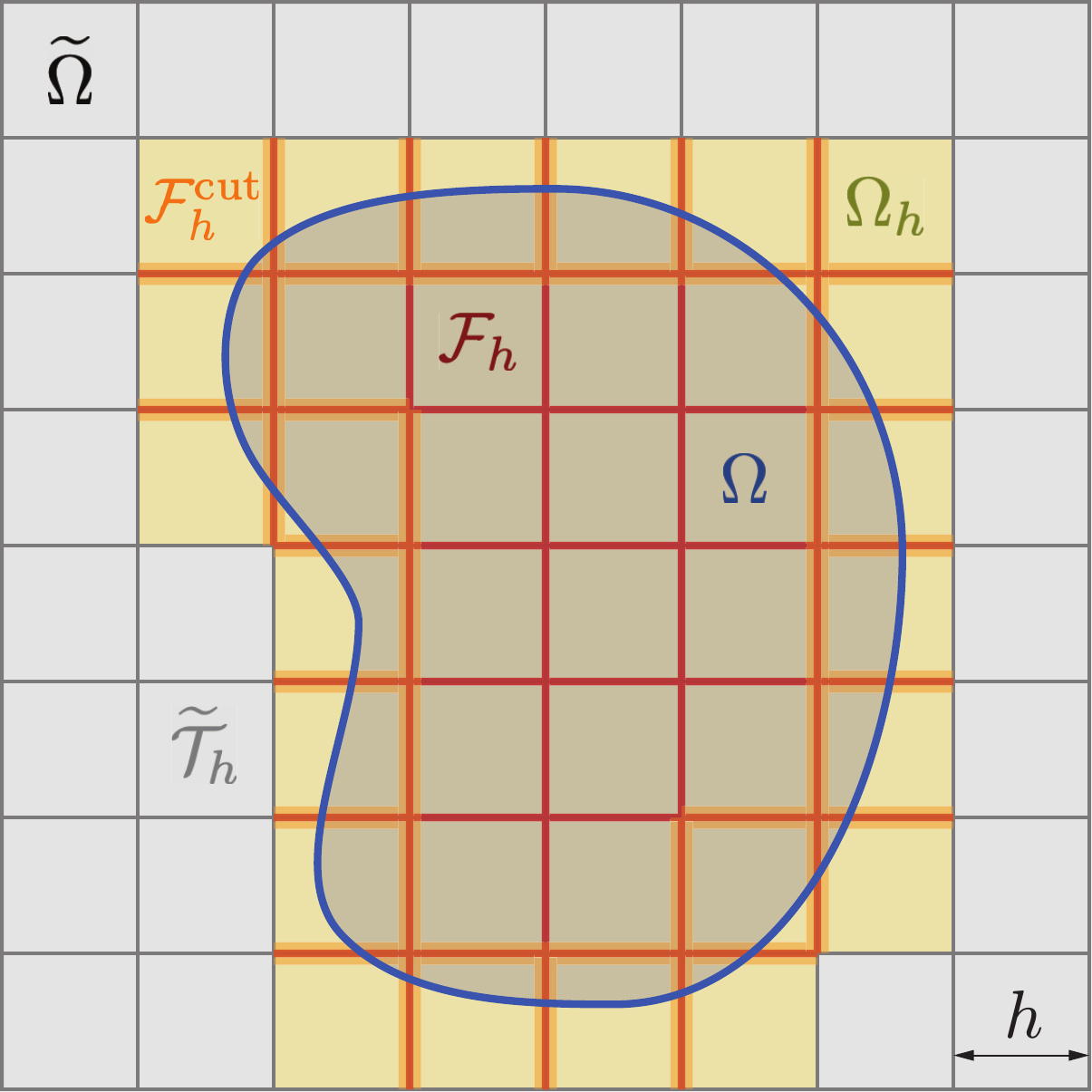}
\caption{CutIGA notation: ambient domain $\widetilde{\Omega}$, background
mesh $\widetilde{\mcT}_h$ with mesh parameter~$h$, physical domain $\Omega$,
active domain $\Omega_h$, skeleton $\mcF_h$, and cut skeleton
$\mcF_h^{\mathrm{cut}}$}
\label{fig:cutiga-notation}
\end{figure}

The active mesh is defined by
\begin{equation}
\mcT_h
:=
\bigl\{
T\in\widetilde{\mcT}_h:T\cap\Omega\neq\emptyset
\bigr\},
\end{equation}
and the corresponding active domain is the geometric realization
$\Omega_h:=\operatorname{int}\bigl(\bigcup_{T\in\mcT_h}\overline{T}\bigr)$; see Figure~\ref{fig:cutiga-notation}.
We distinguish the sets of interior and cut elements according to
\begin{subequations}
\begin{align}
\mcT_h^{\mathrm{int}}
&:=
\bigl\{T\in\mcT_h:T\subset\Omega\bigr\},
\\
\mcT_h^{\mathrm{cut}}
&:=
\bigl\{T\in\mcT_h:T\cap\Gamma\neq\emptyset\bigr\}.
\end{align}
\end{subequations}
The \emph{skeleton} $\mcF_h$ of the active mesh is the collection of its interior faces. The subset
of skeleton faces associated with cut elements is denoted by
\begin{equation*}
\mcF_{h}^{\mathrm{cut}}
:=
\bigl\{
F\in\mcF_h:
F\subset\partial T
\text{ for some }T\in\mcT_h^{\mathrm{cut}}
\bigr\}.
\end{equation*}
The pressure skeleton stabilization will be supported on $\mcF_h$, whereas
the velocity ghost-penalty stabilization will be restricted to
$\mcF_h^{\mathrm{cut}}$. 

To ensure robustness of the CutIGA formulation, we impose the following geometric conditions on the domain and the mesh.

\begin{assumption}[Geometric mesh conditions]
\label{ass:geom}
The domain $\Omega$ and the family of background meshes
$\{\widetilde{\mathcal{T}}_h\}_{h\in\mathcal{H}}$ satisfy the following conditions with constants $\Lambda>0$ and
$L_{\mathrm{path}}\in\mathbb{N}$ independent of $h\in\mcH$:
\begin{enumerate}
  \item[\emph{1.}] \emph{Geometry resolution.}
  For every $T \in \mathcal{T}_h^{\mathrm{cut}}$, there exist a face~$F$ of~$T$
  and a Lipschitz function $g_T:D_T\to\mathbb{R}$ such that $\partial\Omega\cap T$ is the graph of $g_T$
  and $\|\nabla_F g_T\|_{L^\infty(D_T)}\leq\Lambda$,
where $D_T$ is the orthogonal projection of $\partial\Omega\cap T$ onto
  the affine hyperplane containing~$F$ and $\nabla_F$ denotes the tangential gradient in this hyperplane.
  \item[\emph{2.}] \emph{Bounded path length.}
    Every
    $T \in \mathcal{T}_h^{\mathrm{cut}}$ can be connected to some
    $T' \in \mathcal{T}_h^{\mathrm{int}}$ by a face path in $\mathcal{T}_h$
    of length at most $L_{\mathrm{path}}$.
  \item[\emph{3.}] \emph{Local face connectivity.}
    For every $R>0$, there exists $L_R\in\mathbb{N}$, independent
    of~$h\in\mcH$, such that any two elements $T,T'\in\mathcal{T}_h$ with
    $\operatorname{dist}(T,T')\leq Rh$ can be connected by a face path
    in $\mathcal{T}_h$ of length at most~$L_R$.
\end{enumerate}
\end{assumption}
Under the standing assumption that the domain $\Omega$ is bounded Lipschitz (Section~\ref{sec:StokesProblem}), Assumption~\ref{ass:geom}.1 corresponds to a resolution requirement for all members of the mesh family $\{\widetilde{\mcT}_h\}_{h\in\mathcal{H}}$. This condition is required in the trace inequalities below, such as~\eqref{eq:cut-element-trace}. For a bounded Lipschitz domain and a
quasi-uniform, shape-regular mesh family, the uniform interior cone
property ensures Assumption~\ref{ass:geom}.2 for all sufficiently small $h$. This condition is needed for the extension estimates below, e.g.,~\eqref{eq:pressure-extension} and~\eqref{eq:velocity-trace-control}.
Assumption~\ref{ass:geom}.3 excludes configurations in which nearby active
elements are connected only through elements outside the active mesh, e.g.,
across a gap in~$\Omega$ of width comparable to~$h$. For a bounded Lipschitz
domain and a quasi-uniform, shape-regular mesh family, this assumption holds for all
sufficiently small~$h$. This condition is
needed for the jump-controlled interpolation estimate in
Lemma~\ref{lem:jump-controlled-interpolation}.

For $h\in\mathcal H$, let
$\widetilde S_{h,k,r}\subset C^r(\widetilde{\Omega})$ denote the
tensor-product spline space of degree~$k$ in each coordinate direction associated with the background
mesh $\widetilde{\mcT}_h$, where $0\leq r\leq k-1$.
The knot vectors of $\widetilde{S}_{h,k,r}$ have full multiplicity $k+1$ at each component of $\partial\widetilde{\Omega}$.
This ensures that the support of every B-spline basis function of $\widetilde{S}_{h,k,r}$ is contained in the closure of $\widetilde{\Omega}$, and that the restriction of $\widetilde{S}_{h,k,r}$ to every
$T \in \widetilde{\mathcal{T}}_h$ is a complete tensor-product polynomial space $\mathbb{Q}_k(T)$, i.e., the space of polynomials of degree at most~$k$ in each coordinate direction.
The active spline space $S_{h,k,r}
:=
\{
\widetilde v|_{\Omega_h}:\widetilde v\in\widetilde S_{h,k,r}
\}
$
is obtained by restriction to the active domain. We consider
equal-order and equal-regularity velocity and pressure approximation spaces 
$V_h:=[S_{h,k,r}]^d$ and $Q_h:=S_{h,k,r}$. The choice $r=k-1$ corresponds to maximum inter-element regularity.

The bilinear form $b$ in~\eqref{eq:stokes-bilinear-b} does not, in
general, satisfy a discrete inf-sup condition for the equal-order pair
$V_h\times Q_h$. To compensate for this deficiency, we introduce a
skeleton stabilization form. For each face $F\in\mcF_h$, we label the
adjacent elements $T^\pm$ according to an arbitrary but fixed ordering
and orient the unit normal $\bfn_F$ from $T^-$ to $T^+$. The
corresponding traces are denoted by~$(\cdot)^\pm$. The skeleton
stabilization form is given by
\begin{equation}\label{eq:pressure-stabilization}
s_h(p,q)
:=
\tau_s\sum_{F\in\mcF_h}
\sum_{l=r+1}^k h^{2l+1}
\big(
\llbracket\partial_n^l p\rrbracket,
\llbracket\partial_n^l q\rrbracket
\big)_{F},
\end{equation}
where $\tau_s>0$ is a stabilization parameter, $\partial_n^l q:=(\bfn_F\cdot\nabla)^lq$, and the jump across $F$
is defined by
$\left.\llbracket w\rrbracket\right|_F:=w^+-w^-$. For vector- and
tensor-valued functions, the jump is understood componentwise. The
stabilization form is independent of the arbitrary ordering of
$T^\pm$.

Restricting the background spline space to the physical domain may produce
basis functions with arbitrarily small support in $\Omega$, which can lead to instability and ill-conditioning~\cite{Prenter:2023rp}.
To ensure cut-independent stability and conditioning, we introduce the
velocity ghost-penalty form
\begin{equation}\label{eq:velocity-ghost-penalty}
g_h(\bfv,\bfw)
:=
\tau_g\sum_{F\in\mcF_h^{\mathrm{cut}}}
\sum_{l=r}^k h^{2l+1}
\big(
\llbracket\partial_n^l\nabla\bfv\rrbracket,
\llbracket\partial_n^l\nabla\bfw\rrbracket
\big)_{F},
\end{equation}
where $\tau_g>0$ is a stabilization parameter. In contrast to $s_h$, the ghost penalty acts only on the cut skeleton
$\mcF_h^{\mathrm{cut}}$.

Because $S_{h,k,r}\subset C^r$, the pressure
jump $\llbracket\partial_n^l q\rrbracket$ vanishes for
$0\leq l\leq r$, whereas the velocity-gradient jump
$\llbracket\partial_n^l\nabla\bfv\rrbracket$ vanishes for
$0\leq l\leq r-1$. The factors
$h^{2l+1}$ provide $L^2$ scaling for $s_h$ and $H^1$ scaling for $g_h$.
For maximum regularity, $r=k-1$, the pressure sum contains $l=k$, while
the velocity sum contains $l=k-1,k$.

\begin{rem}[Equivalent velocity ghost penalty]
\label{rem:equivalent-velocity-ghost-penalty}
An equivalent, lower-derivative form of $g_h$ is
\begin{equation}\label{eq:velocity-ghost-penalty-normal}
g_h^\bot(\bfv,\bfw)
:=
\tau_g\sum_{F\in\mcF_h^{\mathrm{cut}}}
\sum_{l=r+1}^k h^{2l-1}
\big(
\llbracket\partial_n^l\bfv\rrbracket,
\llbracket\partial_n^l\bfw\rrbracket
\big)_F.
\end{equation}
Indeed, with $\nabla_t=(I-\bfn_F\otimes\bfn_F)\nabla$, the decomposition
$
\llbracket\partial_n^l\nabla\bfv\rrbracket
=\llbracket\partial_n^{l+1}\bfv\rrbracket\otimes\bfn_F
+\nabla_t\llbracket\partial_n^l\bfv\rrbracket
$
and the tangential inverse estimate
$\|\nabla_tw\|_F
\lesssim h^{-1}\|w\|_F$
yield $g_h(\bfv,\bfv)\simeq g_h^\bot(\bfv,\bfv)$, with constants
independent of $h$ and the cut configuration. Thus, $g_h^\bot$ provides
the same control using only normal derivatives of $\bfv$.
\end{rem}

For maximum regularity and $2\mu=1$, $s_h$ and $g_h^\bot$ coincide with the
skeleton and ghost-penalty operators in~\cite{Hoang:2019ty} with
$\tau_s=2\gamma$ and $\tau_g=\tilde{\gamma}/2$, where $\gamma$ and
$\tilde{\gamma}$ denote the stabilization parameters used there. The
constants in the analysis below depend on~$\tau_s$ and~$\tau_g$, but are
independent of~$h$ and the cut configuration.

The essential boundary condition on $\Gamma_{\textsc{D}}$ is imposed
weakly by the symmetric Nitsche method. Combining the consistency and
symmetry terms with the penalty term and the velocity ghost penalty, we
define the bilinear form
\begin{equation}
\begin{aligned}
a_h(\bfv,\bfw)
:={}&a(\bfv,\bfw)+g_h(\bfv,\bfw)
-\big((\nabla_{\mathrm{s}}\bfv)\bfn,\bfw\big)_{\Gamma_{\textsc{D}}}
\\
&-\big(\bfv,(\nabla_{\mathrm{s}}\bfw)\bfn\big)_{\Gamma_{\textsc{D}}}+\beta h^{-1}(\bfv,\bfw)_{\Gamma_{\textsc{D}}},
\end{aligned}
\end{equation}
with a suitable choice of the Nitsche parameter $\beta>0$. To account for the fact that test functions do not vanish at the Dirichlet boundary in the CutIGA formulation, we set
\begin{equation}\label{eq:bh-definition}
b_h(\bfv,q)
:=
b(\bfv,q)+(\bfv,q\bfn)_{\Gamma_{\textsc{D}}},
\end{equation}
with $\bfn$ as the exterior unit normal vector to~$\Gamma$.

Combining the above components, the stabilized CutIGA formulation is given by
\begin{multline}\label{eq:cutiga-formulation}
\text{Find }(\bfu,p)\in{}V_h\times{}Q_h\text{ s.t.}
\\
A_h((\bfu,p),(\bfv,q))=\ell_h(\bfv,q)
\\
\forall (\bfv,q)\in V_h\times Q_h,
\end{multline}
with the discrete bilinear and linear forms
\begin{subequations}\label{eq:discrete-forms}
\begin{align}
A_h((\bfu,p),(\bfv,q))
&:=a_h(\bfu,\bfv)+b_h(\bfv,p)
\notag\\
&\phantom{:=}-b_h(\bfu,q)+s_h(p,q),
\label{eq:discrete-bilinear-form}
\\
\ell_h(\bfv,q)
&:=\ell(\bfv)
-\big(\bfg,(\nabla_{\mathrm{s}}\bfv)\bfn\big)_{\Gamma_{\textsc{D}}}
\notag\\
&\phantom{:=}+\beta h^{-1}(\bfg,\bfv)_{\Gamma_{\textsc{D}}}
-(\bfg,q\bfn)_{\Gamma_{\textsc{D}}}.
\label{eq:discrete-linear-form}
\end{align}   
\end{subequations}
The additional terms in $\ell_h$ relative to $\ell$ result from the weak imposition of
$\bfu=\bfg$: the first two correspond to the adjoint-consistency and
penalty terms of Nitsche's method, while the last accounts for the
nonzero normal trace in $b_h$. 

\section{Stability and well-posedness}
\label{sec:stability}
We next establish the stability of the discrete formulation~\eqref{eq:cutiga-formulation}. The proof
relies on a collection of extension, interpolation, and skeleton-jump
estimates for the active spline spaces. To maintain the focus on the
stability argument, these auxiliary CutIGA results are collected in
Appendix~\ref{app:cutiga-tools} and invoked below as needed.

\subsection{Mesh-dependent norms}
We equip the discrete velocity and pressure spaces with the mesh-dependent
norms
\begin{subequations}
\label{eq:mesh-dependent-norms}
\begin{align}
\tn \bfv \tn^2_h&:= \| \nabla_{\mathrm{s}} \bfv \|^2_\Omega + | \bfv |^2_{g_h} 
+ h \| (\nabla_{\mathrm{s}} \bfv)\bfn \|^2_{\Gamma_{\textsc{D}}}
\notag\\
&\phantom{:=}
+ h^{-1} \| \bfv \|^2_{\Gamma_{\textsc{D}}} 
\label{eq:mesh-dependent-norms-v}
\\
\tn q \tn_{h}^2&:=\| q \|^2_{\Omega}+| q |^2_{s_h}.
\label{eq:mesh-dependent-norms-p}
\end{align}
\end{subequations}
Here,
$|q|_{s_h}^2:=s_h(q,q)$
and
$|\bfv|_{g_h}^2:=g_h(\bfv,\bfv)$
are the pressure skeleton-stabilization and velocity ghost-penalty
seminorms induced by
\eqref{eq:pressure-stabilization} and \eqref{eq:velocity-ghost-penalty}, respectively. We equip the product
space $V_h\times Q_h$ with the norm
$\tn(\bfv,q)\tn_h^2:=\tn\bfv\tn_h^2+\tn q\tn_h^2$. We use the same
notation $\tn\,\cdot\,\tn_h$ for the velocity, pressure, and product
norms; the intended norm is clear from its argument.

\subsection{Continuity}
Continuity of the bilinear form $a_h$ follows from the Cauchy--Schwarz
inequality applied to the volume and ghost-penalty terms and its
weighted form applied to the Nitsche terms:
\begin{align}
|a_h(\bfu,\bfv)|
&\lesssim
\|\nabla_{\mathrm{s}}\bfu\|_\Omega\|\nabla_{\mathrm{s}}\bfv\|_\Omega
+|\bfu|_{g_h}|\bfv|_{g_h}
\notag\\
&\quad
+h^{1/2}\|(\nabla_{\mathrm{s}}\bfu)\bfn\|_{\Gamma_{\textsc{D}}}
 h^{-1/2}\|\bfv\|_{\Gamma_{\textsc{D}}}
\notag\\
&\quad
+h^{-1/2}\|\bfu\|_{\Gamma_{\textsc{D}}}
 h^{1/2}\|(\nabla_{\mathrm{s}}\bfv)\bfn\|_{\Gamma_{\textsc{D}}}
\notag\\
&\quad
+\beta h^{-1/2}\|\bfu\|_{\Gamma_{\textsc{D}}}
 h^{-1/2}\|\bfv\|_{\Gamma_{\textsc{D}}}
\notag\\
&\lesssim \tn\bfu\tn_h\tn\bfv\tn_h.
\label{eq:ah-continuity}
\end{align}

To establish continuity of $b_h$, we first
derive a cut-independent pressure trace estimate. The
stabilization $s_h$ provides the extension estimate
\begin{equation}\label{eq:pressure-extension}
\|q\|_{\Omega_h}^2
\lesssim \|q\|_\Omega^2+|q|_{s_h}^2
=\tn q \tn_h^2.
\end{equation}
The estimate is obtained by propagating control from elements with a
nondegenerate intersection with $\Omega$ to the remaining active
elements along uniformly bounded face paths. This local-to-global
argument is the standard ghost-penalty construction~\cite{Burman:2010km}:
the two-element estimate follows from a Taylor expansion about the common
face~\cite[Lemma~5.1]{Massing:2014}, and is iterated along face
paths~\cite[Proposition~5.1]{Massing:2014}; see
also~\cite[Lemma~6]{Burman:2012nn} and the systematic treatments
in~\cite{Burman:2015fp,Burman:2025sf}.
The standard inverse trace inequality
$\|q\|_{T\cap\Gamma_{\textsc{D}}}\lesssim h^{-1/2}\|q\|_T$ holds for
all $T\in\mcT_h$ such that
$T\cap\Gamma_{\textsc{D}}\neq\emptyset$, with a constant independent
of the cut configuration. Summation over these elements yields
\begin{equation}\label{eq:pressure-trace}
h\|q\|_{\Gamma_{\textsc{D}}}^2
\lesssim
\|q\|_{\Omega_h}^2
\lesssim
\tn q \tn_h^2,
\end{equation}
where the last inequality follows from~\eqref{eq:pressure-extension}.
From
$\|\nabla\cdot\bfv\|_\Omega\lesssim
\|\nabla_{\mathrm{s}}\bfv\|_\Omega$ and \eqref{eq:pressure-trace},
it follows that
\begin{align}
|b_h(\bfv,q)|
&\leq
\|\nabla\cdot\bfv\|_\Omega\|q\|_\Omega
+\|\bfv\|_{\Gamma_{\textsc{D}}}
 \|q\|_{\Gamma_{\textsc{D}}}
\notag\\
&\lesssim \tn\bfv\tn_h\tn q\tn_h.
\label{eq:bh-continuity}
\end{align}

The Cauchy--Schwarz inequality for $s_h$ yields
\begin{equation}\label{eq:sh-continuity}
|s_h(p,q)|\leq |p|_{s_h}|q|_{s_h}
\leq \tn p\tn_h\tn q\tn_h.
\end{equation}
Combining \eqref{eq:ah-continuity}--\eqref{eq:sh-continuity}, we obtain
\begin{equation}\label{eq:Ah-continuity}
|A_h((\bfu,p),(\bfv,q))|
\lesssim
\tn(\bfu,p)\tn_h\tn(\bfv,q)\tn_h.
\end{equation}

To establish continuity of the linear form~$\ell_h$ in
\eqref{eq:discrete-linear-form}, we apply the Cauchy--Schwarz
inequality and its weighted form to obtain
\begin{align}
|\ell_h(\bfv,q)|
&\leq
\|\bff\|_\Omega\|\bfv\|_\Omega
+\|\bfh\|_{\Gamma_{\textsc{N}}}
 \|\bfv\|_{\Gamma_{\textsc{N}}}
\notag\\
&\quad
+h^{-1/2}\|\bfg\|_{\Gamma_{\textsc{D}}}
 h^{1/2}\|(\nabla_{\mathrm{s}}\bfv)\bfn\|_{\Gamma_{\textsc{D}}}
\notag\\
&\quad
+\beta h^{-1/2}\|\bfg\|_{\Gamma_{\textsc{D}}}
 h^{-1/2}\|\bfv\|_{\Gamma_{\textsc{D}}}
\notag\\
&\quad
+h^{-1/2}\|\bfg\|_{\Gamma_{\textsc{D}}}
 h^{1/2}\|q\|_{\Gamma_{\textsc{D}}}.
\end{align}
Because $\Gamma_{\textsc{D}}$ has positive measure, Korn's inequality
with a boundary term states that
\begin{equation}\label{eq:korn-boundary}
\|\bfv\|_{H^1(\Omega)}
\lesssim
\|\nabla_{\mathrm{s}}\bfv\|_\Omega
+\|\bfv\|_{\Gamma_{\textsc{D}}},
\quad \bfv\in \bfH^1(\Omega).
\end{equation}
The boundary term controls the so-called rigid-body motions in the kernel of the
symmetric gradient. Using the definition of $\tn\bfv\tn_h$ and the
trace theorem, \eqref{eq:korn-boundary} therefore implies
\begin{equation}\label{eq:korn-trace}
\|\bfv\|_{H^1(\Omega)}
+\|\bfv\|_{\Gamma_{\textsc{N}}}
\lesssim \tn\bfv\tn_h,
\qquad \bfv\in V_h.
\end{equation}
Using \eqref{eq:pressure-trace} and \eqref{eq:korn-trace}, we conclude
that
\begin{equation}\label{eq:lh-continuity}
|\ell_h(\bfv,q)|
\lesssim
\|(\bff,\bfh,\bfg)\|_h
\tn(\bfv,q)\tn_h,
\end{equation}
with the mesh-dependent data norm
\begin{equation}
\|(\bff,\bfh,\bfg)\|_h^2
:=
\|\bff\|_\Omega^2
+\|\bfh\|_{\Gamma_{\textsc{N}}}^2
+h^{-1}\|\bfg\|_{\Gamma_{\textsc{D}}}^2.
\end{equation}
\begin{rem}
The continuity constant is independent of $h$ and the cut configuration,
but for fixed nonzero $\bfg$, the Nitsche scaling causes the data norm
to grow like $h^{-1/2}$ as $h\to0$. Uniform-in-$h$ estimates require an
$\bfH^1(\Omega)$-stable lifting of
$\bfg\in\bfH^{1/2}(\Gamma_{\textsc{D}})$ or homogeneous Dirichlet data.
\end{rem}

\subsection{Coercivity}
We first establish coercivity of $a_h$. For $\bfv\in V_h$ and any
$\delta>0$, Young's inequality gives
\begin{align*}
a_h(\bfv,\bfv)
&=\|\nabla_{\mathrm{s}}\bfv\|_\Omega^2+|\bfv|_{g_h}^2
-2\big((\nabla_{\mathrm{s}}\bfv)\bfn,\bfv\big)_{\Gamma_{\textsc{D}}}
\\
&\phantom{=}+\beta h^{-1}\|\bfv\|_{\Gamma_{\textsc{D}}}^2
\\
&\geq \|\nabla_{\mathrm{s}}\bfv\|_\Omega^2+|\bfv|_{g_h}^2
-\delta h\|(\nabla_{\mathrm{s}}\bfv)\bfn\|_{\Gamma_{\textsc{D}}}^2
\\
&\phantom{=}+(\beta-\delta^{-1})h^{-1}\|\bfv\|_{\Gamma_{\textsc{D}}}^2.
\end{align*}
Analogously to~\eqref{eq:pressure-trace}, the elementwise inverse-trace
inequality combined with the ghost-penalty extension estimate, applied to the
piecewise-polynomial tensor~$\nabla_{\mathrm{s}}\bfv$, whose jumps are bounded
by those of~$\nabla\bfv$, yields
\begin{equation}\label{eq:velocity-trace-control}
h\|(\nabla_{\mathrm{s}}\bfv)\bfn\|_{\Gamma_{\textsc{D}}}^2
\lesssim \|\nabla_{\mathrm{s}}\bfv\|_\Omega^2+|\bfv|_{g_h}^2.
\end{equation}
Hence, choosing $\delta>0$ sufficiently small and then $\beta$
sufficiently large absorbs the two boundary contributions and leads to
\begin{equation}\label{eq:velocity-coercivity}
a_h(\bfv,\bfv) \gtrsim \tn \bfv \tn_h^2,
\qquad \bfv\in V_h,
\end{equation}
uniformly with respect to the cut configuration.

We next use Verf\"urth's trick~\cite{Verfurth:1984}, as in~\cite{Burman:2006wf}, to establish weak
coercivity of~$b_h$ in~\eqref{eq:bh-definition}.  The Neumann boundary induces a pressure-trace term. To account for the fact that
$\Gamma_{\textsc{N}}$ may intersect the background mesh arbitrarily, this term requires special treatment.
Under mild assumptions on the domain geometry and mesh resolution,
Appendix~\ref{app:neumann-boundary-neighborhood} establishes an estimate
on a discrete Neumann-boundary neighborhood that yields the required
pressure-trace estimate in Lemma~\ref{lem:neumann-pressure-trace} below.
Lemma~\ref{lem:relaxed-infsup} identifies the 
pressure-gradient and trace contributions to the
LBB deficiency, which are controlled by
Lemmas~\ref{lem:lbb-deficiency} and
\ref{lem:neumann-pressure-trace}, respectively. Proposition~\ref{prop:pressure-infsup} presents the stabilized pressure inf-sup condition.

\begin{lemma}[Relaxed inf-sup estimate]\label{lem:relaxed-infsup}
There are constants $c_1,c_2,c_3>0$ such that, for all $q\in Q_h$,
\begin{equation}\label{eq:verfurth-a}
c_1\|q\|_\Omega-c_2h\|\nabla q\|_\Omega
-c_3h^{1/2}\|q\|_{\Gamma_{\textsc{N}}}
\leq
\sup_{\bfv\in V_h}\frac{b_h(\bfv,q)}{\tn\bfv\tn_h}.
\end{equation}
\end{lemma}
\begin{proof}
The standard inf-sup condition for the bilinear form
$b:\bfH^1_{\Gamma_{\textsc{D}},0}(\Omega)\times L^2(\Omega)\to\IR$
implies the existence of a
constant $c_1>0$ such that, for every $q\in Q_h\subset L^2(\Omega)$,
there exists $\bfv\in\bfH^1_{\Gamma_{\textsc{D}},0}(\Omega)$ satisfying
\begin{equation}
c_1\|q\|_\Omega
\leq
\frac{b_h(\bfv,q)}{\|\nabla\bfv\|_\Omega}.
\end{equation}
Let $\pi_h:L^2(\Omega)\to{}S_{h,k,r}$ denote the componentwise interpolant defined in
Appendix~\ref{app:cutiga-tools}. The approximation and stability properties of $\pi_h$, in combination with the cut-element trace inequality 
\begin{equation}\label{eq:cut-element-trace}
\|v\|_{T\cap\Gamma}^2
\lesssim
h^{-1}\|v\|_T^2+h\|\nabla v\|_T^2,
\quad v\in{}H^1(T),\; T\in\mcT_h,
\end{equation}
which holds uniformly in the cut configuration by
Assumption~\ref{ass:geom}.1, give
\begin{subequations}
\begin{align}
\|\bfv-\pi_h\bfv\|_\Omega
+h^{1/2}\|\bfv-\pi_h\bfv\|_{\Gamma_{\textsc{N}}}
&\lesssim h\|\nabla\bfv\|_\Omega,
\\
\tn\pi_h\bfv\tn_h&\lesssim\|\nabla\bfv\|_\Omega.
\end{align}
\end{subequations}
Integration by parts yields
\begin{equation*}
b_h(\pi_h\bfv-\bfv,q)
=(\pi_h\bfv-\bfv,\nabla q)_\Omega
-((\pi_h\bfv-\bfv),q\bfn)_{\Gamma_{\textsc{N}}}.
\end{equation*}
Therefore,
\begin{multline*}
b_h(\pi_h\bfv,q)
=b_h(\bfv,q)+b_h(\pi_h\bfv-\bfv,q)
\\
\geq
\bigl(c_1\|q\|_\Omega-c_2h\|\nabla q\|_\Omega
-c_3h^{1/2}\|q\|_{\Gamma_{\textsc{N}}}\bigr)
\|\nabla\bfv\|_\Omega.
\end{multline*}
If the parenthesized expression is positive,
estimate~\eqref{eq:verfurth-a} follows from
\(\tn\pi_h\bfv\tn_h\lesssim\|\nabla\bfv\|_\Omega\). If it is nonpositive,
the estimate follows from the nonnegativity of the supremum.
\end{proof}

\begin{rem}
The negative terms on the left-hand side of~\eqref{eq:verfurth-a}
represent the {\em LBB deficiency\/}~\cite{Bochev:2006gf}.
\end{rem}

\begin{lemma}[Pressure-trace estimate]
\label{lem:neumann-pressure-trace}
Suppose Assumption~\ref{ass:inward-progressing-patches} in Appendix~\ref{app:neumann-boundary-neighborhood} holds. For every
$\epsilon>0$, there is a constant $C_\epsilon>0$ such that, after reducing the baseline mesh width~$h_0$ if necessary,
\begin{equation}\label{eq:neumann-pressure-trace}
h^{1/2}\|q\|_{\Gamma_{\textsc{N}}}
\leq
\epsilon\|q\|_\Omega+C_\epsilon|q|_{s_h},
\quad q\in Q_h,\:h\in\mathcal H.
\end{equation}
\end{lemma}

\begin{proof}
Let $U_h$ denote the union of the active elements intersecting
$\Gamma_{\textsc{N}}$.  The standard inverse-trace inequality gives
$h^{1/2}\|q\|_{\Gamma_{\textsc{N}}}\lesssim\|q\|_{U_h}$.
By Lemma~\ref{lem:boundary-neighborhood}, for every $\delta>0$ and all
sufficiently small $h$,
\begin{equation}
\label{eq:pressure-boundary-neighborhood}
\|q\|_{U_h}^2
\leq
\delta\|q\|_\Omega^2+C_\delta|q|_{s_h}^2.
\end{equation}
Taking square roots and choosing $\delta$ sufficiently small in terms
of $\epsilon$ proves~\eqref{eq:neumann-pressure-trace}.
\end{proof}

The skeleton stabilization~\eqref{eq:pressure-stabilization} only involves
jumps of normal derivatives. To relate it to jumps of the full pressure
gradient, we introduce the jump seminorm
\begin{equation}\label{eq:jh}
|v|_{j_h}^2
:=
\sum_{F\in\mcF_h}\sum_{l=0}^k
h^{2l+1}\big\|\llbracket\nabla^l v\rrbracket\big\|_F^2.
\end{equation}

\begin{lemma}[Equivalence of jump seminorms]\label{lem:seminormequiv}
For all $q\in Q_h$, it holds that
\begin{equation}
\label{eq:seminormequiv}
h|\nabla q|_{j_h}\simeq|q|_{s_h}.
\end{equation}
\end{lemma}

\begin{proof}
By definition,
\begin{equation}\label{eq:gradient-jump-reindexed}
h^2|\nabla{}q|_{j_h}^2=
\sum_{F\in\mcF_h}\sum_{l=1}^{k+1}
h^{2l+1}\big\|\llbracket\nabla^l q\rrbracket\big\|_F^2.
\end{equation}
Consider an arbitrary face~$F$ and let $\{\bft_1,\ldots,\bft_{d-1}\}$ be an orthonormal basis for its tangent plane, aligned with the edges of~$F$. For a multi-index
$\boldsymbol{\alpha}=(\alpha_1,\ldots,\alpha_{d-1})$, we define
\begin{equation}
\partial_t^{\boldsymbol{\alpha}}
=(\boldsymbol{t}_1\cdot\nabla)^{\alpha_1}
(\boldsymbol{t}_2\cdot\nabla)^{\alpha_2}
\cdots
(\boldsymbol{t}_{d-1}\cdot\nabla)^{\alpha_{d-1}}.
\end{equation}
In the orthonormal frame $\{\bft_1,\ldots,\bft_{d-1},\bfn_F\}$, the jump $\smash[t]{\llbracket\nabla^lq\rrbracket}$
in~\eqref{eq:gradient-jump-reindexed} has components of the form
$\llbracket\partial_t^{\boldsymbol{\alpha}}\partial_n^aq\rrbracket$
with $|\boldsymbol{\alpha}|+a=l$, where
$|\boldsymbol{\alpha}|=\alpha_1+\cdots+\alpha_{d-1}$.
Since tangential differentiation commutes with taking the jump, we have
$\llbracket\partial_t^{\boldsymbol{\alpha}}\partial_n^aq\rrbracket
=
\partial_t^{\boldsymbol{\alpha}}\llbracket\partial_n^aq\rrbracket$.
Applying the standard polynomial inverse inequality on~$F$ and using
$l=|\boldsymbol{\alpha}|+a$, we obtain
\begin{equation}
\label{eq:estimate}
\begin{aligned}
h^{2l+1}\big\|\llbracket\partial_t^{\boldsymbol{\alpha}}\partial_n^aq\rrbracket\big\|_F^2
&\lesssim h^{2l+1-2|\boldsymbol{\alpha}|}
\big\|\llbracket\partial_n^aq\rrbracket\big\|_F^2
\\
&=
h^{2a+1}
\big\|\llbracket\partial_n^aq\rrbracket\big\|_F^2.
\end{aligned}
\end{equation}
In addition, for $q\in Q_h$, the expression
$\partial_t^{\boldsymbol{\alpha}}\llbracket\partial_n^aq\rrbracket$
vanishes if $a\leq r$, if $a>k$, or if $\alpha_i>k$ for some index~$i$. Summing the estimate~\eqref{eq:estimate} over the
components, derivative orders, and faces in~\eqref{eq:gradient-jump-reindexed}
yields $h|\nabla q|_{j_h}\lesssim|q|_{s_h}$.
Conversely, since $\partial_n^aq$ is a component of $\nabla^aq$,
summing the corresponding terms for $a=r+1,\ldots,k$ gives
$|q|_{s_h}\lesssim h|\nabla q|_{j_h}$.
\end{proof}

\begin{lemma}[Control of the LBB deficiency]\label{lem:lbb-deficiency}
There are constants $c_4,c_5,c_6>0$ such that, for all $q\in Q_h$,
\begin{equation}\label{eq:verfurth-b}
c_4h\|\nabla q\|_\Omega-c_5|q|_{s_h}
-c_6h^{1/2}\|q\|_{\Gamma_{\textsc{N}}}
\leq
\sup_{\bfv\in V_h}\frac{b_h(\bfv,q)}{\tn\bfv\tn_h}.
\end{equation}
\end{lemma}

\begin{proof}
Let $\tilde{q}\in\widetilde{S}_{h,k,r}$ be such that $\tilde{q}|_{\Omega_h}=q$,
and set $\bfz:=(\pi_h^*\nabla\tilde{q})|_{\Omega_h}$, with $\pi_h^*$ applied
componentwise, and $\bfw:=h^2\bfz$. Since $(\pi_h^*w)|_{\Omega_h}$ depends only
on $w|_{\Omega_h}$ (see Appendix~\ref{app:cutiga-tools}), $\bfz$ does not
depend on the choice of~$\tilde{q}$.
Combining an active-mesh extension estimate for $\nabla q$, analogous
to~\eqref{eq:pressure-extension}, with the stability
$\|\bfz\|_{\Omega_h}\lesssim\|\nabla q\|_{\Omega_h}$ of the
interpolant, which follows from~\eqref{eq:oswald}, yields
\begin{equation}
\|\bfz\|_{\Omega_h}
\lesssim \|\nabla q\|_\Omega+|\nabla q|_{j_h}.
\end{equation}
The inverse trace inequality and standard polynomial inverse estimates
therefore imply
\begin{subequations}
\begin{align}
h^{1/2}\|\bfz\|_{\Gamma_{\textsc{N}}}
&\lesssim \|\nabla q\|_\Omega+|\nabla q|_{j_h},
\\
\tn\bfw\tn_h
&\lesssim h\bigl(\|\nabla q\|_\Omega
+|\nabla q|_{j_h}\bigr).
\label{eq:gradient-test-stability}
\end{align}  
\end{subequations}
Moreover, the jump-controlled interpolation estimate
\eqref{eq:jump-controlled-interpolation} for piecewise polynomial functions,
established in Appendix~\ref{app:cutiga-tools}, gives
\begin{equation}
\|\bfz-\nabla q\|_\Omega
\lesssim |\nabla q|_{j_h}.
\end{equation}
Using integration by parts, we obtain
\begin{align}
b_h(\bfw,q)
&=h^2\big((\bfz,\nabla q)_\Omega
-(\bfz,q\bfn)_{\Gamma_{\textsc{N}}}\big)
\notag\\
&=h^2\big(\|\nabla{}q\|_{\Omega}^2+(\bfz-\nabla{}q,\nabla{}q)_{\Omega}
-(\bfz,q\bfn)_{\Gamma_{\textsc{N}}}\big)
\notag\\
&\geq h^2\|\nabla q\|_\Omega^2
-C h^2|\nabla q|_{j_h}\|\nabla q\|_\Omega
\notag\\
&\phantom{=}
-C h^{3/2}\bigl(\|\nabla q\|_\Omega
+|\nabla q|_{j_h}\bigr)
\|q\|_{\Gamma_{\textsc{N}}},
\label{eq:bhwq}
\end{align}
for some $C>0$.
From~\eqref{eq:gradient-test-stability} and~\eqref{eq:bhwq} it follows that
\begin{multline}
\bigg(\sup_{\bfv\in V_h}\frac{b_h(\bfv,q)}{\tn\bfv\tn_h}+
C h^{1/2}\|q\|_{\Gamma_{\textsc{N}}}
\bigg)
\bigl(\|\nabla q\|_\Omega
+|\nabla q|_{j_h}\bigr)
\\
+
Ch|\nabla q|_{j_h}\|\nabla q\|_\Omega
\geq h\|\nabla q\|_\Omega^2.
\end{multline}
Young's inequality then yields
\begin{equation}\label{eq:pressure-gradient-young}
h^2\|\nabla q\|_\Omega^2
\lesssim
\bigg(\sup_{\bfv\in V_h}\frac{b_h(\bfv,q)}{\tn\bfv\tn_h}\bigg)^2
+h^2|\nabla q|_{j_h}^2
+h\|q\|_{\Gamma_{\textsc{N}}}^2.
\end{equation}
Inequality~\eqref{eq:verfurth-b} follows from~\eqref{eq:pressure-gradient-young}, the triangle inequality, and Lemma~\ref{lem:seminormequiv}.
\end{proof}

\begin{proposition}[Stabilized pressure inf-sup estimate]
\label{prop:pressure-infsup}
Suppose Assumption~\ref{ass:inward-progressing-patches} holds.  There are
constants $c_7,c_8>0$ such that, for all $q\in Q_h$,
\begin{equation}\label{eq:infsup-small}
c_7 \|q\|_\Omega-c_8|q|_{s_h}
\leq
\sup_{\bfv\in V_h}\frac{b_h(\bfv,q)}{\tn\bfv\tn_h}.
\end{equation}
\end{proposition}
\begin{proof}
Taking a convex combination of
\eqref{eq:verfurth-a} and \eqref{eq:verfurth-b} with weights
$\zeta$ and $1-\zeta$, respectively, where $0\leq\zeta\leq1$, and applying
\eqref{eq:neumann-pressure-trace} yields
\begin{multline*}
\zeta c_1 \|q\|_{\Omega} +((1-\zeta)c_4-\zeta c_2)h\|\nabla q\|_{\Omega}
-(1-\zeta)c_5|q|_{s_h}
\\
-(\zeta c_3+(1-\zeta)c_6)(\epsilon\|q\|_{\Omega}+C_{\epsilon}|q|_{s_h})
\leq
\sup_{\bfv\in V_h}\frac{b_h(\bfv,q)}{\tn\bfv\tn_h}.
\end{multline*}
Choose $0<\zeta<1$ such that
$(1-\zeta)c_4-\zeta c_2\geq0$. Next, choose $\epsilon>0$
sufficiently small that
$c_7:=\zeta c_1-(\zeta c_3+(1-\zeta)c_6)\epsilon>0$,
reducing $h_0$ as required by Lemma~\ref{lem:neumann-pressure-trace}.
Setting
$c_8:=(1-\zeta)c_5+(\zeta c_3+(1-\zeta)c_6)C_\epsilon$
and dropping the nonnegative pressure-gradient term then
gives~\eqref{eq:infsup-small}.
\end{proof}

\begin{corollary}[Pure Dirichlet problem]
\label{cor:pure-dirichlet}
Suppose that $\Gamma_{\textsc{N}}=\emptyset$, and let
$Q_{h,0}$ denote the subspace
of $Q_h$ of functions with vanishing mean over~$\Omega$. There are
constants $c_7,c_8>0$ such that~\eqref{eq:infsup-small} holds for all
$q\in Q_{h,0}$, without Assumption~\ref{ass:inward-progressing-patches}.
\end{corollary}
\begin{proof}
For $\Gamma_{\textsc{N}}=\emptyset$, the standard inf-sup condition holds
for $b:\bfH^1_0(\Omega)\times L^2_0(\Omega)\to\IR$, where $L^2_0(\Omega)$
denotes the subspace of $L^2(\Omega)$ of functions with vanishing mean
over~$\Omega$. The proofs of
Lemmas~\ref{lem:relaxed-infsup} and~\ref{lem:lbb-deficiency} then carry over
for $q\in Q_{h,0}$ without the Neumann-boundary terms, i.e., with
$c_3=c_6=0$, and Lemma~\ref{lem:neumann-pressure-trace} is not required.
Choosing $0<\zeta<1$ such that $(1-\zeta)c_4-\zeta c_2\geq0$ in the proof of
Proposition~\ref{prop:pressure-infsup} yields~\eqref{eq:infsup-small} with
$c_7:=\zeta c_1$ and $c_8:=(1-\zeta)c_5$.
\end{proof}

\subsection{Well-posedness and stability}
We combine the velocity coercivity estimate
\eqref{eq:velocity-coercivity} with the stabilized pressure inf-sup
estimate \eqref{eq:infsup-small} to establish well-posedness and stability of~\eqref{eq:cutiga-formulation}. 

\begin{proposition}[Inf--sup stability and well-posedness]
\label{prop:discrete-stability}
Under Assumption~\ref{ass:inward-progressing-patches}, there is a constant~$\gamma>0$,
independent of the cut configuration, such that for all
$(\bfv,q)\in V_h\times Q_h$,
\begin{equation}
\label{eq:infsup-big}
\gamma\tn(\bfv,q)\tn_h
\leq
\sup_{0\neq(\bfw,r)\in (V_h\times Q_h)}
\frac{A_h((\bfv,q),(\bfw,r))}{\tn(\bfw,r)\tn_h}.
\end{equation}
Consequently, the discrete problem~\eqref{eq:cutiga-formulation} admits a unique solution
$(\bfu,p)\in V_h\times Q_h$, which satisfies the stability estimate
\begin{equation}\label{eq:solution-stability}
\tn(\bfu,p)\tn_h
\leq C_\ell\gamma^{-1}\|(\bff,\bfh,\bfg)\|_h,
\end{equation}
where $C_\ell>0$ denotes the hidden constant in~\eqref{eq:lh-continuity},
which is likewise independent of the cut configuration.
\end{proposition}

\begin{proof}
By \eqref{eq:infsup-small} and Young's inequality, there are constants
$c_p,C_p>0$ such that, for each $q\in Q_h$, there exists $\bfz\in V_h$
with $\tn\bfz\tn_h=\|q\|_\Omega$ satisfying
\begin{equation}\label{eq:pressure-test-bound}
b_h(\bfz,q)
\geq c_p\|q\|_\Omega^2-C_p|q|_{s_h}^2.
\end{equation}
We use the test pair $(\bfw,r)=(\bfv+\sigma\bfz,q)$ with $\sigma>0$. Let $c_a>0$ and $C_a>0$ denote the hidden constants in the coercivity estimate~\eqref{eq:velocity-coercivity} and the continuity estimate~\eqref{eq:ah-continuity}, respectively.
Expanding $A_h$ according to~\eqref{eq:discrete-bilinear-form} and
using~\eqref{eq:pressure-test-bound}, we obtain
\begin{equation}
\begin{aligned}
&A_h((\bfv,q),(\bfv+\sigma\bfz,q))
\\
&=a_h(\bfv,\bfv)+|q|_{s_h}^2
+\sigma{}b_h(\bfz,q)
+\sigma{}a_h(\bfv,\bfz)
\\
&\geq c_a\tn\bfv\tn_h^2+(1-\sigma{}C_p)|q|_{s_h}^2
+\sigma c_p\|q\|_\Omega^2
\\
&\phantom{\geq}
-\sigma C_a\tn\bfv\tn_h\|q\|_\Omega.
\end{aligned}  
\end{equation}
Applying Young's inequality to the mixed term and then choosing
$\sigma>0$ sufficiently small yields
\begin{equation}\label{eq:full-norm-lower-bound}
\tn\bfv\tn_h^2+|q|_{s_h}^2+\|q\|_\Omega^2
\lesssim
A_h((\bfv,q),(\bfv+\sigma\bfz,q)).
\end{equation}
The left-hand side of~\eqref{eq:full-norm-lower-bound} equals
$\tn(\bfv,q)\tn_h^2$. Moreover, the triangle inequality and
$\tn\bfz\tn_h=\|q\|_\Omega$ imply
\begin{equation}
\tn(\bfv+\sigma\bfz,q)\tn_h
\lesssim \tn(\bfv,q)\tn_h.
\end{equation}
Combining these bounds gives~\eqref{eq:infsup-big}. Since $V_h\times Q_h$
is finite dimensional, estimate~\eqref{eq:infsup-big} implies that the discrete problem~\eqref{eq:cutiga-formulation} admits a unique solution. The stability estimate~\eqref{eq:solution-stability} follows by applying~\eqref{eq:infsup-big} to the
solution $(\bfu,p)$ of the CutIGA formulation~\eqref{eq:cutiga-formulation}
and invoking the continuity estimate~\eqref{eq:lh-continuity}.
\end{proof}

\section{A priori error estimates}
\label{sec:a-priori-error-estimates}
The error analysis follows the structure of Strang's first lemma for
conforming but non-consistent approximations~\cite{ErnGuermond2004}, with one
modification: the role of the original bilinear form is taken by a modified
bilinear form~$A_h^*$, which is consistent and continuous on a space that
contains both the exact solution and $V_h\times Q_h$. The form $A_h^*$ is
obtained from $A_h$ in~\eqref{eq:discrete-bilinear-form} by omitting the
stabilization terms, which involve jumps of higher-order derivatives and are
generally not defined for the solution of~\eqref{eq:stokes-weak-formulation}.
The Nitsche terms are retained; to ensure that these are well defined, we
assume that the solution of~\eqref{eq:stokes-weak-formulation} resides in
$\bfH^2(\Omega)\times{}H^1(\Omega)$. The difference $A_h-A_h^*$ then
represents the consistency error, which is controlled by the stabilization
seminorms. We introduce $A_h^*$ in tandem with the corresponding starred norms
$\tn\cdot\tn_h^*$:
\begin{align*}
{\tn\bfv\tn_h^*}^2
&:=\|\nabla_{\mathrm{s}}\bfv\|_\Omega^2
+h\|(\nabla_{\mathrm{s}}\bfv)\bfn\|_{\Gamma_{\textsc{D}}}^2
+h^{-1}\|\bfv\|_{\Gamma_{\textsc{D}}}^2,
\\
{\tn q\tn_h^*}^2
&:=\|q\|_\Omega^2+h\|q\|_{\Gamma_{\textsc{D}}}^2,
\end{align*}
and the corresponding product norm
${\tn(\bfv,q)\tn_h^*}^2
:={\tn\bfv\tn_h^*}^2+{\tn q\tn_h^*}^2$. 
The pressure norm includes the scaled trace contribution $h\|q\|_{\Gamma_{\textsc{D}}}^2$ to control the Dirichlet-boundary pressure trace in the velocity--pressure coupling~$b_h$. The pressure-trace estimate~\eqref{eq:pressure-trace} and the norm definitions imply
\begin{equation}
\label{eq:star-norm-equivalence}
\tn(\bfv,q)\tn_h^*\lesssim\tn(\bfv,q)\tn_h,
\end{equation}
for all $(\bfv,q)\in V_h\times Q_h$.

\subsection{Quasi-optimal approximation}
The modified bilinear form $A_h^*$ is consistent in the sense that every
sufficiently smooth solution $(\widehat{\bfu},\widehat{p})$
of~\eqref{eq:stokes-weak-formulation} satisfies, for all $(\bfv,q)\in V_h\times Q_h$,
\begin{equation}\label{eq:Ah-star-consistency}
A^*_h((\widehat{\bfu},\widehat{p}),(\bfv,q))=\ell_h(\bfv,q).
\end{equation}
Moreover, it is continuous according to
\begin{equation}
\label{eq:Ah-star-continuity}
\big|A^*_h((\bfu,p),(\bfv,q))\big|\lesssim\tn(\bfu,p)\tn^*_h\tn(\bfv,q)\tn_h,
\end{equation}
for all $(\bfu,p)\in{}(\boldsymbol{H}^2(\Omega)\times H^1(\Omega))+(V_h\times Q_h)$ and $(\bfv,q)\in V_h\times Q_h$.

\begin{proposition}[Quasi-optimal approximation]
\label{prop:quasi-optimal-approximation}
Let $(\widehat{\bfu},\widehat{p})\in\boldsymbol{H}^2(\Omega)\times H^1(\Omega)$ and $(\bfu,p)\in V_h\times Q_h$ denote the solutions to~\eqref{eq:stokes-weak-formulation} and~\eqref{eq:cutiga-formulation}, respectively. The following quasi-optimal approximation property holds:
\begin{multline}
\label{eq:quasi-optimal-approximation}
\tn(\widehat{\bfu},\widehat{p})-(\bfu,p)\tn_h^*
\leq
\inf_{(\bfv,q)\in V_h\times Q_h}
\bigg(
  \frac{C_0}{\gamma}\big(
|\bfv|_{g_h}+|q|_{s_h} \big) 
\\
+\bigg(1+\frac{C_0C_1}{\gamma}\bigg)\tn(\widehat{\bfu},\widehat{p})-(\bfv,q)\tn_h^*
\bigg).
\end{multline}
Here, $\gamma>0$ is the stability constant in~\eqref{eq:infsup-big}, and $C_0,C_1>0$ are the hidden constants in~\eqref{eq:star-norm-equivalence} and~\eqref{eq:Ah-star-continuity}, respectively.
\end{proposition}
\begin{proof}
Set $\bfU:=(\bfu,p)$ and $\widehat{\bfU}:=(\widehat{\bfu},\widehat{p})$.
The CutIGA formulation~\eqref{eq:cutiga-formulation}, consistency
relation~\eqref{eq:Ah-star-consistency}, and continuity
estimate~\eqref{eq:Ah-star-continuity} imply
\begin{align*}
&|A_h(\bfU-\bfV,\bfW)|
=
|\ell_h(\bfW)-A_h(\bfV,\bfW)|
\\
&=|A_h^*(\widehat{\bfU}-\bfV,\bfW)+A_h^*(\bfV,\bfW)-A_h(\bfV,\bfW)|
\\
&\leq C_1\tn\widehat{\bfU}-\bfV\tn^*_h\tn\bfW\tn_h
+|A_h^*(\bfV,\bfW)-A_h(\bfV,\bfW)|,
\end{align*}
for all $\bfV,\bfW\in V_h \times Q_h$. 
It then follows from the 
inf-sup stability estimate~\eqref{eq:infsup-big} that
\begin{multline}
\gamma\tn\bfU-\bfV\tn_h
\leq
C_1\tn\widehat{\bfU}-\bfV\tn_h^*
\\
+\sup_{0\neq\bfW\in V_h\times Q_h}
\frac{|A_h^*(\bfV,\bfW)-A_h(\bfV,\bfW)|}{\tn\bfW\tn_h}.
\end{multline}
Because $A_h^*$ and $A_h$ differ only in the stabilization terms, and for all $\bfW:=(\bfw,r)\in V_h \times Q_h$ and $(\bfv,q)\in V_h\times Q_h$
it holds that 
$|g_h(\bfv,\bfw)|
\leq |\bfv|_{g_h}\tn\bfW\tn_h$
and
$|s_h(q,r)|
\leq |q|_{s_h}\tn\bfW\tn_h$ by virtue of the Cauchy--Schwarz inequality, the supremum term on the right-hand side is bounded by
$|\bfv|_{g_h}+|q|_{s_h}$. Estimate~\eqref{eq:quasi-optimal-approximation}
then follows by combining the preceding bound with the triangle inequality and~\eqref{eq:star-norm-equivalence},
\begin{equation}
\tn\widehat{\bfU}-\bfU\tn_h^*
\leq
\tn\widehat{\bfU}-\bfV\tn_h^*
+C_0
\tn\bfU-\bfV\tn_h,
\end{equation}
and taking the infimum over~$\bfV\in V_h\times Q_h$.
\end{proof}

\subsection{Optimal-rate convergence}
The quasi-optimal approximation property of the CutIGA formulation according to Proposition~\ref{prop:quasi-optimal-approximation} in combination with the interpolation estimate~\eqref{eq:spline-interpolation-estimate} in Appendix~\ref{app:cutiga-tools} implies optimal-rate convergence of~\eqref{eq:cutiga-formulation}, provided that the stabilization terms 
in~\eqref{eq:quasi-optimal-approximation}
can be controlled.

\begin{lemma}[Stabilization seminorms of interpolants]
\label{lem:stabilization-seminorms-interpolants}
Consider the interpolation operator $\pi^*_h:L^2(\widetilde{\Omega})\to{}\widetilde{S}_{h,k,r}$ introduced in Appendix~\ref{app:cutiga-tools}. For all $\bfv\in\bfH^s(\widetilde{\Omega})$ and
$q\in{}H^t(\widetilde{\Omega})$ with $1\leq s,t \leq k+1$ it holds that
\begin{subequations}
\begin{align}
|\pi^*_h\bfv|_{g_h}&\lesssim h^{s-1}|\bfv|_{H^s(\widetilde{\Omega})},
\label{eq:jump-control-v}
\\
|\pi^*_hq|_{s_h}&\lesssim h^t|q|_{H^t(\widetilde{\Omega})}.
\label{eq:jump-control-p}
\end{align}
\end{subequations}
\end{lemma}
\begin{proof}
Consider a face $F=T^+\cap T^-\in\mcF_h$ and set
$\omega_F:=\omega_{T^+}\cup\omega_{T^-}$. Since $T^+$ and $T^-$ share the
face~$F$, $\omega_F$ is an axis-aligned box with edge lengths $\simeq h$, and
the Bramble--Hilbert lemma provides a polynomial
$p\in\mathbb{P}_{t-1}(\IR^d)\subset\mathbb{Q}_k(\IR^d)$ such that
\begin{equation}
\|q-p\|_{\omega_F}\lesssim h^t|q|_{H^t(\omega_F)}.
\end{equation}
Because $p|_{\widetilde{\Omega}}\in \widetilde{S}_{h,k,r}$, the invariance
property of $\pi^*_h$ yields $\pi^*_hp=p$. Moreover, the smoothness of $p$
implies $\llbracket\partial_n^lp\rrbracket=0$ for all $l\geq{}0$. The inverse
trace inequality on $T^+$ and $T^-$, the local stability
estimate~\eqref{eq:oswald}, and the inclusion
$\omega_{T^\pm}\cap\Omega_h\subset\omega_F$ then give
\begin{multline*}
h^{2l+1}
\bigl\|\llbracket\partial_n^l\pi_h^*q\rrbracket\bigr\|_F^2
=
h^{2l+1}
\bigl\|\llbracket\partial_n^l\pi_h^*(q-p)\rrbracket\bigr\|_F^2
\\
\lesssim
\|\pi_h^*(q-p)\|_{T^+\cup T^-}^2
\lesssim
\|q-p\|_{\omega_F}^2
\lesssim
h^{2t}|q|_{H^t(\omega_F)}^2.
\end{multline*}
Summing over $l=r+1,\ldots,k$ and all faces $F\in\mcF_h$, and using the
bounded overlap of the patches $\omega_F$, we obtain~\eqref{eq:jump-control-p}.
The estimate~\eqref{eq:jump-control-v} follows in a similar manner.
\end{proof}

When the comparison functions in~\eqref{eq:quasi-optimal-approximation} are chosen as interpolants of extensions of the exact solution,
Lemma~\ref{lem:stabilization-seminorms-interpolants} provides optimal-rate control of the stabilization terms. Together with the standard interpolation estimate~\eqref{eq:spline-interpolation-estimate} for the remaining terms, this yields optimal convergence rates.

\begin{proposition}[Optimal-rate convergence]
\label{prop:optimal-rate-convergence}
Let $(\widehat{\bfu},\widehat{p})\in\bfH^{s}(\Omega)\times{}H^t(\Omega)$ and $(\bfu,p)\in{}V_h\times{}Q_h$ denote the solutions to~\eqref{eq:stokes-weak-formulation} and~\eqref{eq:cutiga-formulation}, respectively, with
$2\leq s \leq k+1$ and $1\leq t \leq k+1$. It holds that
\begin{equation}\label{eq:optimal-rate-convergence}
\tn(\widehat{\bfu},\widehat{p})-(\bfu,p)\tn_h^*\lesssim
h^{s-1}\|\widehat{\bfu}\|_{\bfH^s(\Omega)}
+
h^t\|\widehat{p}\|_{H^t(\Omega)}.
\end{equation}
\end{proposition}
\begin{proof}
Let $\widehat{\bfu}{}^E:=E\widehat{\bfu}$ and $\widehat{p}{}^E:=E\widehat{p}$
denote the extensions of the exact solution by the extension operator~$E$
introduced in Appendix~\ref{app:cutiga-tools}, with $E$ applied
componentwise to~$\widehat{\bfu}$. The parenthesized expression on the
right-hand side of~\eqref{eq:quasi-optimal-approximation} can be estimated
by choosing
$(\bfv,q)=(\pi_h\widehat{\bfu},\pi_h\widehat{p})
=(\pi_h^*\widehat{\bfu}{}^E,\pi_h^*\widehat{p}{}^E)|_{\Omega_h}$.
Expanding the norm $\tn\cdot\tn_h^*$ yields
\begin{multline}
 {\bigtn\big(\widehat{\bfu},\widehat{p}\big)-\big(\pi^*_h\widehat{\bfu}{}^E,\pi^*_h\widehat{p}{}^E\big)\bigtn_h^*}^2
 \\
 =
 \big\|\nabla_{\mathrm{s}}\big(\widehat{\bfu}{}^E-\pi^*_h\widehat{\bfu}{}^E\big)\big\|_{\Omega}^2
+
 h\big\|\nabla_{\mathrm{s}}\big(\widehat{\bfu}{}^E-\pi^*_h\widehat{\bfu}{}^E\big)\bfn\big\|_{\Gamma_{\textsc{D}}}^2
 \\
+
h^{-1}\big\|\widehat{\bfu}{}^E-\pi^*_h\widehat{\bfu}{}^E\big\|_{\Gamma_{\textsc{D}}}^2
+
\big\|\widehat{p}{}^E-\pi^*_h\widehat{p}{}^E\big\|_{\Omega}^2
\\
+
h\big\|\widehat{p}{}^E-\pi^*_h\widehat{p}{}^E\big\|_{\Gamma_{\textsc{D}}}^2.
\end{multline}
Since $|\nabla_{\mathrm{s}}\bfv|\leq|\nabla\bfv|$ for all
$\bfv\in\bfH^1(\Omega)$, the interpolation estimate of
Lemma~\ref{lem:spline-interpolation-estimate} with $m=1$ implies that the
first term on the right-hand side is bounded by
$\lesssim h^{2s-2}\|\widehat{\bfu}\|_{\bfH^s(\Omega)}^2$. To estimate the
second term, we select an element $T\in\mcT_h$ such that
$T\cap\Gamma_{\textsc{D}}\neq\emptyset$. The cut-element trace
inequality~\eqref{eq:cut-element-trace} implies
\begin{multline}
\label{eq:interpolation-gradient-boundary}
h\|\nabla_{\mathrm{s}}\big(\widehat{\bfu}{}^E-\pi^*_h\widehat{\bfu}{}^E\big)\bfn\big\|_{T\cap\Gamma_{\textsc{D}}}^2
\lesssim
\big\|\nabla\big(\widehat{\bfu}{}^E-\pi^*_h\widehat{\bfu}{}^E\big)\big\|_T^2
\\
+
h^2\big\|\nabla^2\big(\widehat{\bfu}{}^E-\pi^*_h\widehat{\bfu}{}^E\big)\big\|_T^2.
\end{multline}
Summing~\eqref{eq:interpolation-gradient-boundary} over all active
elements that intersect the Dirichlet boundary and applying
Lemma~\ref{lem:spline-interpolation-estimate} with $m=1,2$ conveys that the
second term is bounded by
$\lesssim h^{2s-2}\|\widehat{\bfu}\|_{\bfH^s(\Omega)}^2$. Applying the same
trace argument to the third term, now using
Lemma~\ref{lem:spline-interpolation-estimate} with $m=0,1$, yields the same
bound. Lemma~\ref{lem:spline-interpolation-estimate} with $m=0$ implies that
the fourth term is bounded by
$\lesssim h^{2t}\|\widehat{p}\|_{H^t(\Omega)}^2$. Finally, considering the
fifth term, the cut-element trace inequality~\eqref{eq:cut-element-trace}
on the active elements intersecting the Dirichlet boundary and Lemma~\ref{lem:spline-interpolation-estimate} with
$m=0,1$ give
\begin{multline*}
h\|\widehat{p}{}^E-\pi^*_h\widehat{p}{}^E\|_{\Gamma_{\textsc{D}}}^2
\lesssim
\sum_{\substack{T\in\mcT_h\\T\cap\Gamma_{\textsc{D}}\neq\emptyset}}
\Big(\|\widehat{p}{}^E-\pi^*_h\widehat{p}{}^E\|_{T}^2
\\
+h^2\|\nabla(\widehat{p}{}^E-\pi^*_h\widehat{p}{}^E)\|_{T}^2\Big)
\lesssim h^{2t}\|\widehat{p}\|_{H^t(\Omega)}^2.
\end{multline*}
Lemma~\ref{lem:stabilization-seminorms-interpolants} and the stability of
the extension operator provide the bounds
$|\pi_h^*\widehat{\bfu}{}^E|_{g_h}\lesssim
h^{s-1}\|\widehat{\bfu}\|_{\bfH^s(\Omega)}$ and
$|\pi_h^*\widehat{p}{}^E|_{s_h}\lesssim h^t\|\widehat{p}\|_{H^t(\Omega)}$
on the stabilization seminorms.
\end{proof}

\begin{rem}[Optimality in the product norm]
Proposition~\ref{prop:optimal-rate-convergence} implies
$\tn(\widehat{\bfu},\widehat{p})-(\bfu,p)\tn_h^*
\lesssim h^k$ for $s=t=k+1$.
This rate is optimal for the solution pair in the product norm,
because the velocity error is measured in an $H^1$-type norm and its
best-approximation rate in the degree-$k$ space $V_h$ is $h^k$. The resulting estimate for the pressure component separately is, however, one order below its best $L^2$-approximation rate: for
$\widehat{p}\in H^{k+1}(\Omega)$, the optimal approximation rate in $Q_h$ is of order
$h^{k+1}$. 
\end{rem}

\begin{rem}[Pure Dirichlet problem]
\label{rem:pure-dirichlet}
For $\Gamma_{\textsc{N}}=\emptyset$, the data must satisfy the compatibility
condition $\int_\Gamma\bfg\cdot\bfn=0$, and the exact pressure is normalized
such that $\widehat{p}\in L^2_0(\Omega)$. Replacing $Q_h$ by $Q_{h,0}$,
Propositions~\ref{prop:discrete-stability},
\ref{prop:quasi-optimal-approximation},
and~\ref{prop:optimal-rate-convergence} carry over, with
Corollary~\ref{cor:pure-dirichlet} in place of
Proposition~\ref{prop:pressure-infsup}. In the proof of
Proposition~\ref{prop:optimal-rate-convergence}, the pressure comparison
function $\pi_h\widehat{p}$ is replaced by its mean-free part; because
$\widehat{p}\in L^2_0(\Omega)$, the mean of $\pi_h\widehat{p}$ is bounded by
$\|\widehat{p}-\pi_h\widehat{p}\|_\Omega$, and constants do not contribute
to~$|\cdot|_{s_h}$, so that the convergence rates are unaffected.
\end{rem}

\section{Condition number}
\label{sec:condition-number}
We conclude the analysis by establishing cut-independence of the condition number of
the matrix associated with the CutIGA formulation. 

Let $\mcB_h$ be the standard B-spline basis of $S_{h,k,r}$.
For $\varphi\in \mcB_h$ and $T\in\mcT_h$, set
\begin{subequations}
\label{eq:ThphiBhT}
\begin{align}
\mcT_{h,\varphi} &= \{ T \in \mcT_h : T \subset \operatorname{supp}(\varphi) \},
\\
\mcB_{h,T} &= \{ \varphi\in\mcB_h : T \in\mcT_{h,\varphi}\}.
\end{align}
\end{subequations}
For any element $T\in\mcT_h$ and coefficient vector
$(a_\varphi)_{\varphi\in\mcB_{h,T}}$, scaling to the reference element and
the local stability of the B-spline basis yield
\begin{equation}
\label{eq:b-spline-stability}
\Big\|\sum_{\varphi\in\mcB_{h,T}}a_{\varphi}\varphi\Big\|_T^2
\simeq
h^d
\sum_{\varphi\in\mcB_{h,T}}a_{\varphi}^2.
\end{equation}
Considering $(\bfv,q)\in V_h\times Q_h$, we collect the
coefficients of $(v_1,\ldots,v_d,q)$ relative to the basis~$\mcB_h$ in a vector
$\mathbf{x}\in\mathbb R^{N_h}$ with $N_h=(d+1)\,\#\mcB_h$.
The local stability property~\eqref{eq:b-spline-stability} and the uniformly bounded overlap of
basis-function supports imply
\begin{equation}
\label{eq:coefficient-l2-equivalence}
h^{d/2}\|\mathbf{
x}\|_{\ell^2}
\simeq
\bigl(\|\bfv\|_{\Omega_h}^2+\|q\|_{\Omega_h}^2\bigr)^{1/2},
\end{equation}
with $\|\cdot\|_{\ell^2}$ the Euclidean coefficient norm.
The hidden constants in~\eqref{eq:coefficient-l2-equivalence} are independent of~$h$ and
of the cut configuration.

The pressure-extension estimate~\eqref{eq:pressure-extension}, the
corresponding velocity ghost-penalty extension estimate, and Korn's
inequality with the Nitsche boundary term give
\begin{equation}
\label{eq:l2-control-by-discrete-norm}
\|\bfv\|_{\Omega_h}+\|q\|_{\Omega_h}
\lesssim
\tn(\bfv,q)\tn_h.
\end{equation}
Conversely, standard polynomial inverse and inverse-trace estimates, together
with the definitions of $g_h$ and $s_h$, yield
\begin{equation}
\label{eq:discrete-norm-inverse-estimate}
\tn(\bfv,q)\tn_h
\lesssim
h^{-1}\bigl(\|\bfv\|_{\Omega_h}^2+\|q\|_{\Omega_h}^2\bigr)^{1/2}.
\end{equation}
Combining~\eqref{eq:coefficient-l2-equivalence}--
\eqref{eq:discrete-norm-inverse-estimate}, we obtain
\begin{equation}
\label{eq:coefficient-discrete-norm-equivalence}
h^{d/2}\|\mathbf{x}\|_{\ell^2}
\lesssim
\tn(\bfv,q)\tn_h
\lesssim
h^{d/2-1}\|\mathbf{x}\|_{\ell^2}.
\end{equation}

\begin{proposition}[Cut-independent condition-number estimate]
\label{prop:condition-number}
Let $\mathbf{A}_h\in\mathbb R^{N_h\times N_h}$ denote the system matrix of~\eqref{eq:cutiga-formulation} defined by
\begin{equation}
\mathbf{y}^{T}\mathbf{A}_h\mathbf{x}
:=
A_h((\bfv,q),(\bfw,r)),
\end{equation}
where $\mathbf{x}$ and $\mathbf{y}$ are the coefficient vectors of
$(\bfv,q)$ and $(\bfw,r)$, respectively.
Under the assumptions of Proposition~\ref{prop:discrete-stability}, the
matrix $\mathbf{A}_h$ is nonsingular and satisfies
\begin{equation}
\label{eq:condition-number}
\kappa_2(\mathbf{A}_h)
:=
\big\|\mathbf{A}_h\big\|_{\ell^2}\big\|\mathbf{A}_h^{-1}\big\|_{\ell^2}
\lesssim h^{-2},
\end{equation}
with a constant independent of the cut configuration.
\end{proposition}

\begin{proof}
The continuity estimate~\eqref{eq:Ah-continuity} and the upper bound
in~\eqref{eq:coefficient-discrete-norm-equivalence} imply
\begin{multline*}
\|\mathbf{A}_h\|_{\ell^2}
:=\sup_{0\neq\mathbf{x}\in\IR^{N_h}}\frac{\|\mathbf{A}_h\mathbf{x}\|_{\ell^2}}{\|\mathbf{x}\|_{\ell^2}}
\\
=
\sup_{0\neq\mathbf{x}\in\IR^{N_h}}\sup_{0\neq\mathbf{y}\in\IR^{N_h}}\frac{|\mathbf{y}^T\mathbf{A}_h\mathbf{x}|}{\|\mathbf{x}\|_{\ell^2}\|\mathbf{y}\|_{\ell^2}}
\lesssim h^{d-2}.
\end{multline*}
On the other hand, 
\begin{multline*}
\big\|\mathbf{A}_h^{-1}\big\|_{\ell^2}
:=
\sup_{0\neq\mathbf{x}\in\IR^{N_h}}\frac{\|\mathbf{A}_h^{-1}\mathbf{x}\|_{\ell^2}}{\|\mathbf{x}\|_{\ell^2}}
\\
=\sup_{0\neq\mathbf{x}\in\IR^{N_h}}\frac{\|\mathbf{x}\|_{\ell^2}}{\|\mathbf{A}_h\mathbf{x}\|_{\ell^2}}
=
\bigg(\inf_{0\neq\mathbf{x}\in\IR^{N_h}}\frac{\|\mathbf{A}_h\mathbf{x}\|_{\ell^2}}{\|\mathbf{x}\|_{\ell^2}}\bigg)^{-1}
\\
=\bigg(\inf_{\rule{0pt}{8pt}0\neq\mathbf{x}\in\IR^{N_h}}\sup_{0\neq\mathbf{y}\in\IR^{N_h}}\frac{|\mathbf{y}^T\mathbf{A}_h\mathbf{x}|}{\|\mathbf{x}\|_{\ell^2}\|\mathbf{y}\|_{\ell^2}}\bigg)^{-1}
\lesssim h^{-d}.
\end{multline*}
The final inequality follows from the lower bound
in~\eqref{eq:coefficient-discrete-norm-equivalence} 
and the inf-sup
estimate~\eqref{eq:infsup-big}.
Combining the bounds for $\|\mathbf{A}_h\|_{\ell^2}$ 
and $\|\mathbf{A}^{-1}_h\|_{\ell^2}$ proves~\eqref{eq:condition-number}.
\end{proof}

\begin{rem}[Conditioning for the pure Dirichlet problem]
\label{rem:pure-dirichlet-conditioning}
The condition-number estimate in Proposition~\ref{prop:condition-number}
does not carry over directly to the pure Dirichlet case of
Corollary~\ref{cor:pure-dirichlet}, because $Q_{h,0}$ does not admit a local
B-spline basis; the conditioning then depends on how the mean-value
constraint is imposed, e.g., by a Lagrange multiplier.
\end{rem}

\section{Numerical experiments}
\label{sec:numerical-experiments}

\begin{figure*}[t]
  \centering
  \includegraphics{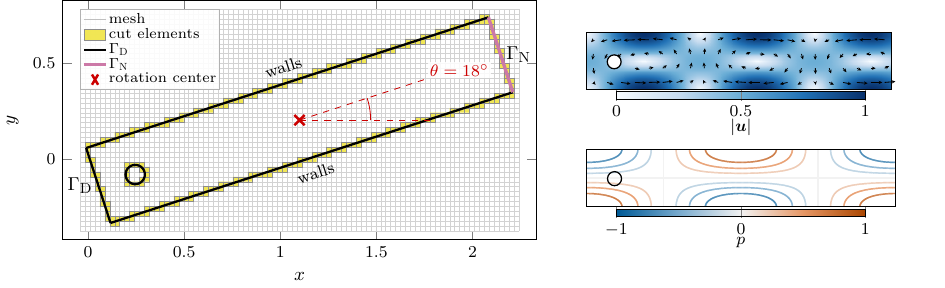}
  \vspace{-14.5pt}
  \caption{Setup, rotation about the channel center (left), and the manufactured solution in the frame of the channel (right)}
  \label{fig:setup}
\end{figure*}

\begin{figure}[t]
  \centering
  \includegraphics{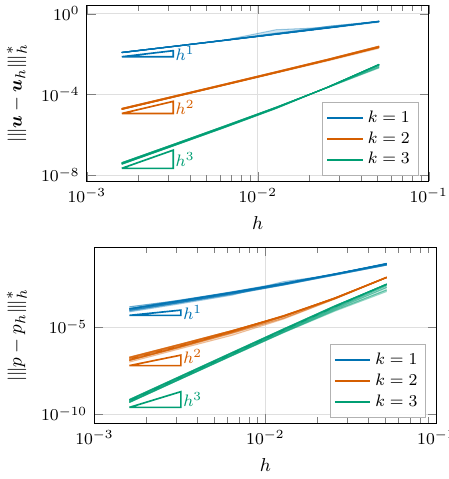}
  \vspace{-13.6pt}
  \caption{Convergence rates of the error $\tn\cdot\tn^*_h$ against the mesh size $h$ for the velocity (top) and pressure (bottom). The thin lines show the error for each angle $\theta$, the dashed line shows the median over all angles, and the shaded area shows the min-max range over all angles. The triangles indicate the expected convergence rates of $h^k$}
  \label{fig:conv}
\end{figure}

\begin{figure}[t]
  \centering
  \includegraphics{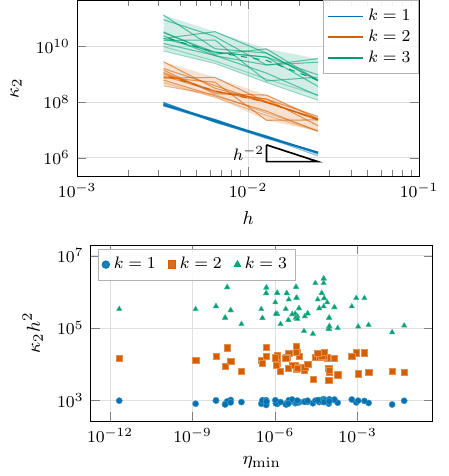}
  \vspace{-12.1pt}
  \caption{Top: condition number $\kappa_2$ against the mesh size $h$, where the thin lines show the condition number for each angle $\theta$, the dashed line shows the median over all angles, and the shaded area shows the min-max range over all angles. Bottom: normalized condition number $\kappa_2h^2$ versus the smallest cut-element volume fraction $\eta_{\min}$}
  \label{fig:channel:kappa}
\end{figure}

This section presents a numerical test to illustrate the convergence rates of
Proposition~\ref{prop:optimal-rate-convergence} and the condition-number
estimate of Proposition~\ref{prop:condition-number}. The test case is based
on~\cite{Hoang:2019ty,Schafer:1996} and illustrated in Figure~\ref{fig:setup}.
We consider a channel of dimensions $L=2.2$, $H=0.41$ containing a cylinder of radius
$R=0.05$ centered at $(0.2,0.2)$. As in~\cite{Prenter:2017rw}, the channel is
rotated about its center $(1.1,0.205)$ with respect to a fixed background mesh,
with angles $\theta\in\{0^\circ,4.5^\circ,9^\circ,\dots,85.5^\circ\}$, giving
20 different cut configurations.
The background mesh has an element size of $h_l=0.41/(8\cdot2^l)$, with refinement levels $l=0,\dots,5$. Velocity and pressure are approximated by maximum-regularity tensor-product B-splines ($r=k-1$) of degree $k=1,2,3$. The Nitsche parameter $\beta=6(k+1)^2$ and the ghost-penalty parameter
$\tilde\gamma=10^{-(k+1)}$ are taken from~\cite{Hoang:2019ty}. For the
skeleton-penalty parameter we use $\gamma=0.1$ for all~$k$, selected by a sweep
over $\gamma\in[10^{-2},10^{-1}]$; this value lies in the range
$\gamma\in(5\cdot10^{-4},5\cdot10^{-1})$ for which the accuracy was found to be
insensitive to~$\gamma$ in~\cite{Hoang:2019ty}.
In the computations, the velocity ghost penalty~$g_h$
in~\eqref{eq:cutiga-formulation} is replaced by the equivalent
form~$g_h^\bot$ in~\eqref{eq:velocity-ghost-penalty-normal}; see
Remark~\ref{rem:equivalent-velocity-ghost-penalty}. The above parameters then
correspond to $\tau_s=0.2$ and $\tau_g=10^{-(k+1)}/2$
in~\eqref{eq:pressure-stabilization} and~\eqref{eq:velocity-ghost-penalty-normal}.
On this domain, we solve the Stokes problem \eqref{eq:stokes} with the manufactured solution:
\begin{equation*}
\begin{gathered}
\psi = b^{-1}\sin(a\eta_1)\sin(b\eta_2),\quad
\bfu = (\partial_{\eta_2}\psi,-\partial_{\eta_1}\psi),\\
p = \cos(c\eta_1)\cos(b\eta_2),\quad
\bff = \mu(a^2+b^2)\bfu+\nabla p,
\end{gathered}
\end{equation*}
with $\boldsymbol\eta\in[0,2.2]\times[0,0.41]$ the coordinates in the frame of the channel, $a=\pi/0.82$, $b=\pi/0.41$, $c=\pi/1.1$ and $\mu=\frac{1}{2}$ as in Section~\ref{sec:StokesProblem}. The exact solution is imposed on the inflow, walls and cylinder using Nitsche's method, and the exact traction is imposed on the outflow. Figure~\ref{fig:setup} (right) shows the manufactured solution. For the integration over cut elements, we use the Duffy--Gauss quadrature rule~\cite{Duffy:1982} for elements cut by the channel walls and the adaptive Algoim rule~\cite{Saye:2015} for elements cut by the cylinder. 

Figure~\ref{fig:conv} displays the velocity and pressure errors in the
norm~$\tn\cdot\tn^*_h$ for $k=1,2,3$ and all considered angles. The figure
conveys that the errors are insensitive to the orientation of the channel with
respect to the background mesh. The velocity errors converge with the optimal rate
$h^k$, in accordance with Proposition~\ref{prop:optimal-rate-convergence}.
The pressure errors converge at least with this rate; over the considered range
of mesh sizes, the observed pressure rates are in fact higher.

Figure~\ref{fig:channel:kappa} presents the condition number
$\kappa_2(\mathbf A_h)$ of the system matrix in~\eqref{eq:condition-number}
on the meshes with $l=1,\dots,4$. Since multiplying the pressure rows of
$\mathbf A_h$ by $-1$ yields a symmetric matrix with the same singular values,
$\kappa_2(\mathbf A_h)$ is computed as the ratio of the largest and smallest
eigenvalues of this matrix in magnitude. For all degrees, the condition number
grows approximately as $h^{-2}$, in accordance with
Proposition~\ref{prop:condition-number}. The spread over the angles increases
with~$k$ but remains bounded. The bottom panel conveys that whereas the smallest cut-element volume fraction
$\eta_{\min}:=\min_{T}|T\cap\Omega|/|T|$~\cite{Prenter:2017rw} ranges over ten
orders of magnitude, the normalized condition number $\kappa_2h^2$ varies by
at most a factor of about~35 for each degree. This confirms that the stabilization removes the
dependence of the conditioning on small cut elements.

\section{Conclusions}
\label{sec:conclusions}
We have presented a stability, error and conditioning analysis of
skeleton-stabilized CutIGA for the Stokes problem with equal-order spline
approximations of the velocity and pressure. For sufficiently fine meshes, the
formulation is inf-sup stable, converges at the optimal rate, and yields a
system matrix whose condition number scales as~$h^{-2}$, with constants
independent of the cut configuration. The analysis applies to
maximum-regularity splines as well as to reduced-regularity spaces down to
standard $C^0$ elements, and covers Neumann boundaries that intersect the mesh
arbitrarily. The numerical experiments confirm the predicted convergence rates
and the cut-independent conditioning.

Several questions remain open. The constants in the analysis depend on the
polynomial degree and on the stabilization parameters, and the experiments
show that the variation of the condition number over the cut configurations
increases with the degree; a degree-explicit analysis could guide the choice of
the stabilization parameters. Furthermore, it may be possible to relax the
geometric condition in Assumption~\ref{ass:inward-progressing-patches}; see
Remark~\ref{rem:inward-patches}.

\bmhead{Acknowledgements}
EHvB\ and MGL\ gratefully acknowledge the hospitality of Institut
Mittag-Leffler in Djursholm, Sweden, during the research program
\emph{Interfaces and Unfitted Discretization Methods} in the fall of 2025,
where parts of this work were discussed.

\appendix
\section{Neumann-boundary neighborhood estimate}
\label{app:neumann-boundary-neighborhood}
This appendix establishes estimate~\eqref{eq:pressure-boundary-neighborhood}.
The proof proceeds in three steps. First, for each active element $T$
intersecting $\Gamma_{\textsc{N}}$, we construct a face-connected patch
extending into the interior of the domain. Second, the skeleton
stabilization restricted to this patch controls the deviation of
$q\in Q_h$ from a single tensor-product polynomial on the patch.
Finally, polynomial mass spreading bounds the norm of this patch
polynomial on $T$ by its norm on suitably selected interior elements.
Combining these estimates and summing over all elements intersecting
$\Gamma_{\textsc{N}}$ yields~\eqref{eq:pressure-boundary-neighborhood}.
The geometric and mesh-resolution assumptions arise from the
construction of the patches.

Consider $q\in{}Q_h$ and two neighboring elements $T^+,T^-\in\mcT_h$ with common face~$F$. The ghost-penalty
estimate~\cite[Lemma~5.1]{Massing:2014} gives
\begin{equation}\label{eq:two-element-pressure-extension}
\|q\|_{T^+}^2
\lesssim
\|q\|_{T^-}^2
+
\sum_{l=r+1}^k h^{2l+1}
\|\llbracket\partial_n^l q\rrbracket\|_F^2.
\end{equation}
Let $\mcP$ be a face-connected element patch, i.e.,
$\mcP\subset\mcT_h$ and, for every $T,T'\in\mcP$, there exist an
integer $N\geq 0$ and a sequence $\{T_j\}_{j=0}^N\subset\mcP$ such
that $T_0=T$, $T_N=T'$, and
$\operatorname{meas}_{d-1}(\partial T_{j-1}\cap\partial T_j)>0$ for
$j=1,\ldots,N$. We denote by $\mcF_{\mcP}$ the collection of interior faces of~$\mcP$:
\begin{multline*}
\mcF_{\mcP}
:=
\big\{\partial T^+\cap\partial{T^-}
  : T^+,T^-\in\mcP,\ T^+\neq T^-,
  \\
  \operatorname{meas}_{d-1}
  (\partial{T^+}\cap\partial{T^-})>0\big\}.
\end{multline*}
The restriction of the skeleton
stabilization~\eqref{eq:pressure-stabilization} to~$\mcP$ is given by:
\begin{equation}\label{eq:local-pressure-seminorm}
|q|_{s_h,\mcP}^2
:=
\tau_s\sum_{F\in\mcF_{\mcP}}
\sum_{l=r+1}^k h^{2l+1}
\big\|\llbracket\partial_n^l q\rrbracket\big\|_F^2.
\end{equation}
Let $\mathbb Q_k(\mathcal P)$ be the space of restrictions to the
geometric realization of~$\mathcal P$ of tensor-product polynomials of
coordinate degree at most~$k$. Iterating~\eqref{eq:two-element-pressure-extension} over the
faces~$F\in\mcF_{\mcP}$ yields
\begin{equation}\label{eq:patch-pressure-reconstruction}
\inf_{p\in\mathbb Q_k(\mathcal P)}
\sum_{T\in\mathcal P}\|q-p\|_T^2
\leq C_{\mathcal P}|q|_{s_h,\mathcal P}^2.
\end{equation}
The constant $C_{\mcP}$ generally depends on the face-adjacency depth of the patch.

Next, we consider polynomial mass spreading. To elucidate the mechanism, we regard a univariate polynomial $p$ of fixed degree $k$ on the unit interval. By virtue of the standard norm equivalence for polynomials~\cite{Schwab1998}, we have
\begin{equation}
\|p\|_{L^\infty(0,1)}^2
\leq C_k\|p\|_{L^2(0,1)}^2.
\end{equation}
Moreover, for any subinterval $\omega\subset(0,1)$, it holds that
$\|p\|_{L^2(\omega)}^2
\leq\operatorname{meas}(\omega)\|p\|_{L^{\infty}(0,1)}^2$
and therefore
\begin{equation}\label{eq:polynomial-mass-spreading}
\|p\|_{L^2(0,\delta)}^2
\leq
C_k\delta\|p\|_{L^2(0,1)}^2,
\end{equation}
for every $\delta\in(0,1]$. The bound~\eqref{eq:polynomial-mass-spreading}
conveys that a fixed-degree polynomial cannot concentrate a fixed fraction of its $L^2$ mass in an arbitrarily small subinterval. Conversely, for any subinterval $(0,\delta)$, the mass of a polynomial in this subinterval spreads to the entire unit interval in accordance with~\eqref{eq:polynomial-mass-spreading}. 

The polynomial spreading argument extends directly to
tensor-product polynomials on multidimensional columns of elements.
For every integer $m\geq 0$, let
$T_j:=D_h\times((0,h)+jh)$
with
$D_h:=(0,h)^{d-1}$ and
$j=0,\ldots,m$. Consider a tensor-product polynomial $p$ of coordinate degree at
most~$k$ on the entire column. For every $x'\in D_h$, the function
$t\mapsto p(x',t)$ is a univariate polynomial of degree at most~$k$.
Applying~\eqref{eq:polynomial-mass-spreading} after rescaling the
interval $(0,(m+1)h)$ to $(0,1)$, and then integrating with respect
to $x'$, yields
\begin{align}\label{eq:column-polynomial-mass-spreading}
\|p\|_{T_0}^2
&=
\int_{D_h}\int_0^h |p(x',t)|^2\,dt\,dx'
\notag\\
&\leq
\frac{C_k}{m+1}
\sum_{j=0}^m
\int_{D_h}\int_{jh}^{(j+1)h}|p(x',t)|^2\,dt\,dx'
\notag\\
&=
\frac{C_k}{m+1}\sum_{j=0}^m\|p\|_{T_j}^2.
\end{align}
For sufficiently large~$m$, the first term in the sum on the right-hand side in~\eqref{eq:column-polynomial-mass-spreading} can be absorbed in the left-hand side and the constant $C_k$ can be adjusted to obtain the bound
\begin{equation}\label{eq:column-polynomial-mass-spreading-absorbed}
\|p\|_{T_0}^2
\leq
\frac{C_k}{m}\sum_{j=1}^m\|p\|^2_{T_j}.
\end{equation}
The bound~\eqref{eq:column-polynomial-mass-spreading-absorbed} conveys
that the squared norm of $p$ on the root element~$T_0$ is controlled
by the sum of its squared norms on the remaining elements, with a
coefficient that can be made arbitrarily small by increasing the
number of elements in the column.

The bound~\eqref{eq:pressure-boundary-neighborhood} follows
from~\eqref{eq:patch-pressure-reconstruction}
and~\eqref{eq:column-polynomial-mass-spreading-absorbed} if every element
intersecting~$\Gamma_{\textsc{N}}$ is the root~$T_0$ of a column
$T_0,\ldots,T_m$ as in~\eqref{eq:column-polynomial-mass-spreading} whose
remaining elements $T_1,\ldots,T_m$ are contained in~$\Omega$. Since these
columns are aligned with the mesh, this assumption imposes restrictive
conditions on both the geometry of~$\Omega$ and the resolution of~$\mcT_h$.
This requirement can be relaxed by considering sufficiently slender
patches progressing into the domain that need not be aligned with the mesh.

\begin{assumption}[Slender inward patches]
\label{ass:inward-progressing-patches}
There exists $m_0\in\mathbb N$ such that, for every integer $m\geq m_0$, there exists $h_m>0$ such that, for
$h\in\mathcal H$ with $h\leq h_m$, every element $T_0\in\mcT_h$
intersecting $\Gamma_{\textsc{N}}$ is associated with 
\begin{itemize}
\item an inward-pointing unit vector $\bfd(T_0)$;
\item a connected patch $\mcP_m(T_0)$;
\item an ordered collection
\begin{equation}
\mcI_m(T_0)=\{S_1,\ldots,S_m\}
\subset\mcT_h^{\mathrm{int}},
\end{equation}
of pairwise distinct elements,
\end{itemize}
such that, for positive constants $c_{\mathrm{in}}$, $C_{\mathrm{in}}$,
$C_{\mathrm{lat}}$, $L_*$, $C_{\mathrm{size}}$, and $C_{\mathrm{spr}}$,
independent of $m$, $h$, $T_0$, and the cut configuration, the following
properties hold:
\begin{enumerate}
\item The barycenters $\bfx_{S_j}$ and $\bfx_{T_0}$ satisfy
\begin{subequations}\label{eq:monotone-inward-column}
\begin{align}
c_{\mathrm{in}}jh
\leq
\big|P_{\bfd(T_0)}(\bfx_{S_j}-\bfx_{T_0})\big|
&\leq C_{\mathrm{in}}jh,
\label{eq:monotone-inward-progress}
\\
\big|P^{\bot}_{\bfd(T_0)}(\bfx_{S_j}-\bfx_{T_0})\big|
&\leq C_{\mathrm{lat}}h,
\label{eq:bounded-lateral-drift}
\end{align}
\end{subequations}
for $j=1,\ldots,m$, with $P_{\bfd(T_0)}=\bfd(T_0)\otimes\bfd(T_0)$ and $P^{\bot}_{\bfd(T_0)}=I-P_{\bfd(T_0)}$  the projection onto $\bfd(T_0)$ and its orthogonal complement, respectively.

\item The element $T_0$ is connected to $S_1$, and $S_j$ is
connected to $S_{j+1}$, by face paths containing at most $L_*$ faces. The intermediate elements may be cut. All elements in these paths belong to $\mathcal P_m(T_0)$, and
\begin{equation}
\#\mathcal P_m(T_0)\leq C_{\mathrm{size}}m.
\end{equation}

\item Every $p\in\mathbb Q_k(\mathcal P_m(T_0))$ satisfies
\begin{equation}\label{eq:patch-polynomial-spreading}
\|p\|_{T_0}^2
\leq
\frac{C_{\mathrm{spr}}}{m}
\sum_{S\in\mathcal I_m}\|p\|_S^2.
\end{equation}

\item The collections $\mathcal I_m(T_0)$ have overlap bounded
independently of $m$, $h$, and the cut configuration.  The face sets
$\mathcal F_{\mathcal P_m(T_0)}$ have overlap bounded by a constant
that may depend on $m$ but not on $h$ or the cut configuration.
\end{enumerate}
\end{assumption}

Conditions~\eqref{eq:monotone-inward-column} arrange the selected elements in $\mathcal I_m(T_0)$ in an $O(h)$-wide tubular
neighborhood of the line
$\bfx_{T_0}+\operatorname{span}(\bfd(T_0))$. Since the connecting face paths have uniformly bounded length, the same is true, up to an
enlargement of the tube width, for the full patch
$\mathcal P_m(T_0)$. Moreover,
\begin{equation}
c_{\mathrm{in}}mh
\leq |\bfx_{S_m}-\bfx_{T_0}|\lesssim mh.
\end{equation}
Thus, the patch $\mathcal P_m(T_0)$ reaches a distance comparable to
$mh$ from the cut element $T_0$. The progression and lateral-drift
bounds in~\eqref{eq:monotone-inward-column}, imposed for every $j$,
control folding and clustering of the selected sequence.

\begin{rem}
\label{rem:inward-patches}
For a bounded Lipschitz domain, we expect
Assumption~\ref{ass:inward-progressing-patches} to hold for all sufficiently
small~$h$: items~1, 2, and~4 follow from the uniform cone condition by
selecting the elements along the cone axes, and item~3 follows from item~1 by
a Remez-type inequality for polynomials on convex
bodies~\cite[Theorem~3 and Remark~3]{BrudnyiGanzburg1973}. A rigorous proof is, however, technical
and beyond the scope of the present work, and we therefore retain
Assumption~\ref{ass:inward-progressing-patches} as a separate assumption.
\end{rem}

\begin{lemma}[Neumann-boundary neighborhood estimate]
\label{lem:boundary-neighborhood}
Suppose Assumption~\ref{ass:inward-progressing-patches} holds. Let
$$
\mcU_h:=\bigl\{T\in\mcT_h:
T\cap\Gamma_{\textsc{N}}\neq\emptyset\bigr\},
$$
and $U_h:=\bigcup_{T\in\mcU_h}T$.
For every $\delta>0$, there is a constant $C_\delta>0$ such that,
after, if necessary, reducing the baseline mesh width~$h_0$,
\begin{equation}\label{eq:small-neumann-boundary-neighborhood}
\|q\|_{U_h}^2
\leq
\delta\|q\|_\Omega^2+C_\delta|q|_{s_h}^2,
\qquad q\in Q_h,\quad h\in\mathcal H.
\end{equation}
\end{lemma}

\begin{proof}
Fix an element $T_0\in\mcU_h$. By
\eqref{eq:patch-pressure-reconstruction}, there exists a single
polynomial $p_m\in\mathbb Q_k(\mathcal P_m(T_0))$ such that
\begin{equation}\label{eq:local-patch-reconstruction}
\sum_{T\in\mathcal P_m(T_0)}\|q-p_m\|_T^2
\leq C_m|q|_{s_h,\mathcal P_m(T_0)}^2.
\end{equation}
From the triangle inequality, the polynomial mass spreading property
\eqref{eq:patch-polynomial-spreading}, and
\eqref{eq:local-patch-reconstruction}, it follows that
\begin{align}
\|q\|_{T_0}^2
&\lesssim
\|q-p_m\|_{T_0}^2+\|p_m\|_{T_0}^2
\notag\\
&\lesssim
C_m|q|_{s_h,\mathcal P_m(T_0)}^2
+\frac{1}{m}\sum_{T\in\mathcal I_m(T_0)}\|p_m\|_T^2
\notag\\
&\lesssim
\frac{1}{m}\sum_{T\in\mathcal I_m(T_0)}\|q\|_T^2
+C_m|q|_{s_h,\mathcal P_m(T_0)}^2.
\label{eq:local-boundary-cell-control}
\end{align}
Summing~\eqref{eq:local-boundary-cell-control} over all
$T_0\in\mcU_h$ and invoking the overlap bounds in
Assumption~\ref{ass:inward-progressing-patches}, we obtain
\begin{equation}\label{eq:m-layer-boundary-neighborhood}
\|q\|_{U_h}^2
\leq
\frac{C}{m}\|q\|_\Omega^2+C_m|q|_{s_h}^2.
\end{equation}
Finally, we choose $m$ so that $C/m\leq\delta$, subsequently reduce $h_0$ so
that $h_0\leq h_m$ and absorb the constant $C_m$ into $C_\delta$ to obtain~\eqref{eq:small-neumann-boundary-neighborhood}.
\end{proof}

\section{CutIGA analysis tools}
\label{app:cutiga-tools}

This appendix collects the interpolation and skeleton-jump estimates
used in the stability and error analysis.

Let $\widetilde{\mcB}_h$ denote the standard B-spline basis of $\widetilde{S}_{h,k,r}$. Define $\widetilde{\mcT}_{h,\varphi}$ and $\widetilde{\mcB}_{h,T}$ as in~\eqref{eq:ThphiBhT}, but with the full background mesh $\widetilde{\mcT}_h$ in place of the active mesh $\mcT_h$. To define the interpolation operator $\pi_h : L^2(\Omega) \to S_{h,k,r}$, we
first extend $v\in L^2(\Omega)$ to $v^E := Ev \in L^2(\widetilde{\Omega})$, where
$E$ denotes Stein's universal extension operator~\cite{Stein1970}, which is bounded from
$H^s(\Omega)$ into $H^s(\widetilde{\Omega})$ for all $s\geq 0$.
For each $T\in\widetilde{\mcT}_h$, let
$P_{T,k}:L^2(T)\rightarrow\mathbb{Q}_k(T)$ be the $L^2(T)$ projection. For
arbitrary $w\in L^2(\widetilde{\Omega})$, the expansion of $P_{T,k}(w|_T)$ in the
basis $\widetilde{\mcB}_{h,T}$ is
\begin{equation}
\label{eq:PTk}
P_{T,k}(w|_T) = \sum_{\varphi\in\widetilde{\mcB}_{h,T}}\hatw_{T,\varphi}\,\varphi|_{T}.
\end{equation}
Thus, for each $\varphi\in\widetilde{\mcB}_h$ and
$T\in\widetilde{\mcT}_{h,\varphi}$, we obtain a local coefficient
$\hatw_{T,\varphi}$. We define the interpolation operator
$\pi^*_h:L^2(\widetilde{\Omega})\to\widetilde{S}_{h,k,r}$ by the map
\begin{equation}
\label{eq:pi*h}
\pi_h^*:w\mapsto\sum_{\varphi\in\widetilde{\mcB}_h}
\langle\hatw_{T,\varphi}\rangle\,\varphi,
\end{equation}
where $\hatw_{T,\varphi}$ are the coefficients in~\eqref{eq:PTk} and
$\langle\cdot\rangle$ is the convex combination
\begin{equation}
\label{eq:<hatW>}
\langle\hatw_{T,\varphi}\rangle
:=
\sum_{T\in\widetilde{\mcT}_{h,\varphi}}\kappa_{T,\varphi}\,\hatw_{T,\varphi}.
\end{equation}
For each $\varphi\in\widetilde{\mcB}_h$, the weights $\kappa_{T,\varphi}\in[0,1]$ form a partition of unity and are otherwise
arbitrary, except that for active basis functions we require
\begin{equation}
\label{eq:kappaT}
\kappa_{T,\varphi} = 0
\qquad
\forall\,\varphi\in\mcB_h,\;
T\in\widetilde{\mcT}_{h,\varphi}\setminus\mcT_h .
\end{equation}
Condition~\eqref{eq:kappaT} is admissible, because
$\mcT_{h,\varphi}=\widetilde{\mcT}_{h,\varphi}\cap\mcT_h\neq\emptyset$ for every
$\varphi\in\mcB_h$. It implies that the coefficient
$\langle\hatw_{T,\varphi}\rangle$ in~\eqref{eq:pi*h} of every active basis
function $\varphi\in\mcB_h$ depends only on the local coefficients
$\hatw_{T,\varphi}$ associated with active elements $T\in\mcT_{h,\varphi}$.
Consequently, $(\pi_h^*w)|_{\Omega_h}$ depends only on $w|_{\Omega_h}$.
Finally, we define the interpolant $\pi_h v := \pi^*_h(v^E)|_{\Omega_h}$.

To establish that the spline space $\widetilde{S}_{h,k,r}$ is invariant under the interpolation operator $\pi_h^*$, we consider arbitrary $v\in{}\widetilde{S}_{h,k,r}$ and~$T\in\widetilde{\mcT}_h$. Evidently, it holds that $P_{T,k} (v|_T) = v|_T$. Therefore, expanding $v$ in the spline basis and restricting to~$T$, we have
\begin{equation*}
 \sum_{\varphi \in \widetilde{\mcB}_h} \hatv_\varphi \varphi|_T = v|_T = P_{T,k}(v|_T) = \sum_{\varphi \in \widetilde{\mcB}_{h,T}} \widehat{v}_{T,\varphi} \varphi|_T.
\end{equation*}
Since the expansion coefficients are unique, we conclude that $\widehat{v}_{T,\varphi} = \hatv_\varphi$, for all $T \in \widetilde{\mcT}_{h,\varphi}$, 
and therefore $\langle\widehat{v}_{T,\varphi}\rangle=\hatv_{\varphi}$. Hence, $\pi^*_h (v) = v$ 
for all~$v \in \widetilde{S}_{h,k,r}$. Note that this argument holds for any choice of convex weights, so the
restriction~\eqref{eq:kappaT} does not affect the invariance property.

The interpolant $\pi_h^*$ is locally $L^2$-stable on active elements: for
$w\in L^2(\widetilde{\Omega})$ and $T\in\mcT_h$,
\begin{equation}
\label{eq:oswald}
\|\pi_h^* w\|_{T} \lesssim \|w\|_{\omega_T\cap\Omega_h},
\end{equation}
where $\omega_T=\bigcup\{\operatorname{supp}\varphi:\varphi\in\widetilde{\mcB}_{h,T}\}$.
In particular, $\|\pi_h^* w\|_{T}\lesssim\|w\|_{\omega_T}$. For $T\in\mcT_h$, it
holds that $\widetilde{\mcB}_{h,T}=\mcB_{h,T}$. Hence, by~\eqref{eq:kappaT},
the coefficients of all basis functions contributing to $(\pi_h^*w)|_T$
involve only local coefficients $\hatw_{T',\varphi}$ on active elements
$T'\subset\omega_T\cap\Omega_h$. Estimate~\eqref{eq:oswald} then follows from
the bound $|\hatw_{T',\varphi}|\lesssim h^{-d/2}\|w\|_{T'}$, which is a
consequence of the $L^2(T')$-stability of $P_{T',k}$ and the local stability
of the B-spline basis~\eqref{eq:b-spline-stability}, together with
$\|\varphi\|_T\lesssim h^{d/2}$ and $\#\mcB_{h,T}\lesssim 1$.

\begin{rem} The estimate~\eqref{eq:oswald} may be viewed as a generalization of the corresponding
result for Oswald interpolation~\cite{Oswald1993}. 
\end{rem}

The invariance property and the local stability~\eqref{eq:oswald} yield the
following interpolation estimate for~$\pi_h$ on the active mesh.

\begin{lemma}[Interpolation estimate]
\label{lem:spline-interpolation-estimate}
For all $v\in H^s(\Omega)$ and $0\le m\le s\le k+1$, it holds that
\begin{equation}
\label{eq:spline-interpolation-estimate}
\sum_{T\in\mcT_h}|v^E-\pi_h v|_{H^m(T)}^2
\lesssim h^{2(s-m)}\|v\|_{H^s(\Omega)}^2,
\end{equation}
where $v^E=Ev$.
\end{lemma}

\begin{proof}
Consider $T\in\mcT_h$. By the tensor-product structure
of~$\widetilde{S}_{h,k,r}$ and the quasi-uniformity of~$\widetilde{\mcT}_h$,
the patch $\omega_T$ is an axis-aligned box containing~$T$ with edge lengths
$\simeq h$. Hence, $\omega_T$ is convex with bounded aspect
ratio, and the Bramble--Hilbert lemma provides a polynomial
$p\in\mathbb{P}_{s-1}(\IR^d)\subset\mathbb{Q}_k(\IR^d)$ such that, for $0\le l\le s$,
\begin{equation}
\label{eq:BH}
|v^E-p|_{H^l(\omega_T)}\lesssim h^{s-l}\,|v^E|_{H^s(\omega_T)},
\end{equation}
with constants independent of $h$, $T$, and the cut configuration. Since
$p|_{\widetilde{\Omega}}\in\widetilde{S}_{h,k,r}$, the invariance property gives
$\pi_h^*p=p$. The triangle inequality therefore yields
\begin{equation}
\label{eq:tri2}
\begin{aligned}
|v^E-\pi_h v|_{H^m(T)}
&\le
|v^E-p|_{H^m(T)}
\\
&\phantom{\le}{}+|\pi_h^*(v^E-p)|_{H^m(T)}.
\end{aligned}
\end{equation}
By~\eqref{eq:BH} with $l=m$ and $T\subset\omega_T$, the first term is bounded
by $h^{s-m}|v^E|_{H^s(\omega_T)}$. For the second term, the inverse inequality
on~$\mathbb{Q}_k(T)$, the stability estimate~\eqref{eq:oswald}, and
\eqref{eq:BH} with $l=0$ give
\begin{equation*}
\begin{aligned}
|\pi_h^*(v^E-p)|_{H^m(T)}
&\lesssim h^{-m}\|\pi_h^*(v^E-p)\|_T
\\
&\lesssim h^{-m}\|v^E-p\|_{\omega_T\cap\Omega_h}
\\
&\le h^{-m}\|v^E-p\|_{\omega_T}
\\
&\lesssim h^{s-m}|v^E|_{H^s(\omega_T)}.
\end{aligned}
\end{equation*}
Note that the Bramble--Hilbert lemma is applied on the convex
patch~$\omega_T$, on which $v^E$ is defined, and not on
$\omega_T\cap\Omega_h$, which may not satisfy its geometric requirements.
Condition~\eqref{eq:kappaT} enters only through the stability
estimate~\eqref{eq:oswald}. Summing over $T\in\mcT_h$ and using
the bounded overlap of the patches~$\omega_T$, we obtain
\begin{equation*}
\sum_{T\in\mcT_h}|v^E-\pi_h v|_{H^m(T)}^2
\lesssim h^{2(s-m)}|v^E|_{H^s(\widetilde{\Omega})}^2.
\end{equation*}
The stability of the extension operator,
$\|v^E\|_{H^s(\widetilde{\Omega})}\lesssim\|v\|_{H^s(\Omega)}$, then
yields~\eqref{eq:spline-interpolation-estimate}.
\end{proof}

The proof of Lemma~\ref{lem:lbb-deficiency} uses an interpolant of
$\nabla q$, with $q\in S_{h,k,r}$, and requires the interpolation error to be controlled by $|\nabla q|_{j_h}$.

\begin{lemma}[Jump-controlled interpolation estimate]
\label{lem:jump-controlled-interpolation}
Let $\widetilde{S}_{h,k,-1}:=\bigoplus_{T\in\widetilde{\mcT}_h}\mathbb{Q}_k(T)$
denote the space of discontinuous piecewise tensor-product polynomials of
degree at most~$k$ in each coordinate direction on the background mesh. For
all $v\in\widetilde{S}_{h,k,-1}$ it holds that
\begin{equation}
\label{eq:jump-controlled-interpolation}
\|v-\pi^*_h v\|_{\Omega_h}\lesssim|v|_{j_h},
\end{equation}
where the jump seminorm $|\cdot|_{j_h}$ is defined in~\eqref{eq:jh}.
\end{lemma}

\begin{proof}
Consider $v\in\widetilde{S}_{h,k,-1}$ and $T\in\mcT_h$, and let
$p\in\mathbb{Q}_k(\IR^d)$ be such that $p|_T=v|_T$. The invariance property
of $\pi^*_h$ yields $\pi^*_h p=p$. Because $v-p$ vanishes on~$T$, the local
stability estimate~\eqref{eq:oswald} implies
\begin{equation*}
\|v-\pi^*_hv\|_T=\|\pi^*_h(v-p)\|_T\lesssim\|v-p\|_{\omega_T\cap\Omega_h}.
\end{equation*}
The set $\omega_T\cap\Omega_h$ is the union of the active elements
$T'\subset\omega_T$, and each such element satisfies
$\operatorname{dist}(T,T')\le kh$. By Assumption~\ref{ass:geom}.3 with $R=k$,
$T'$ can be connected to~$T$ by a face path $T=T_0,T_1,\ldots,T_N=T'$
in~$\mcT_h$ with $N\le L_k$. For neighboring elements $T_{j-1},T_j$ with
common face $F_j\in\mcF_h$, a standard ghost-penalty estimate,
cf.~\eqref{eq:two-element-pressure-extension}, gives for all
$w\in\mathbb{Q}_k(T_{j-1})\oplus\mathbb{Q}_k(T_j)$
\begin{equation*}
\|w\|^2_{T_j}\lesssim\|w\|^2_{T_{j-1}}
+\sum_{l=0}^k h^{2l+1}\|\llbracket\partial_n^l w\rrbracket\|^2_{F_j}.
\end{equation*}
Applying this estimate to $w=v-p$, which vanishes on~$T_0$ and satisfies
$\llbracket\partial_n^l(v-p)\rrbracket=\llbracket\partial_n^l v\rrbracket$,
and iterating along the path yields
\begin{equation*}
\|v-p\|^2_{T'}\lesssim\sum_{F\in\mcF_{T,T'}}\sum_{l=0}^k
h^{2l+1}\|\llbracket\nabla^l v\rrbracket\|_F^2,
\end{equation*}
where $\mcF_{T,T'}\subset\mcF_h$ denotes the faces of the path. Summing over
the uniformly bounded number of active elements $T'\subset\omega_T$, and
subsequently over $T\in\mcT_h$, yields~\eqref{eq:jump-controlled-interpolation}.
Indeed, all faces involved belong to~$\mcF_h$ and lie within a distance
$\lesssim(k+L_k)h$ of~$T$, so that each face is counted a uniformly bounded
number of times.
\end{proof}

\begin{rem} Note that the $h$-scaling in $j_h$ is chosen such that we have the inverse inequality
\begin{equation}
|v|_{j_h} \lesssim \|v\|_{\Omega_h}, \qquad v \in \widetilde{S}_{h,k,-1}.
\end{equation} 
\end{rem}

\def\doiurl#1{\unskip}
\bibliography{BibFile}

\end{document}